\documentclass[twoside]{article}
\usepackage[latin1]{inputenc}
\usepackage[T1]{fontenc}
\usepackage{parskip}
\usepackage{amsmath,amsfonts,amssymb,amsxtra}
\usepackage{latexsym}
\usepackage{theorem}
\usepackage{graphicx,psfrag}
\usepackage{verbatim}
\usepackage{hyperref}

\newtheorem{theo}{Theorem}
\newtheorem{prop}{Proposition}
\newtheorem{lemma}{Lemma}
\newtheorem{cor}{Corollary}

{\theoremstyle{plain}                       
\theorembodyfont{\rmfamily}
}

{\theoremstyle{plain}                       
\theorembodyfont{\rmfamily}
\newtheorem{Example}{Example}}

{\theoremstyle{plain}                       
\theorembodyfont{\rmfamily}
\newtheorem{Definition}{Definition}}

{\theoremstyle{plain}                       
\theorembodyfont{\rmfamily}
}

\def\A{\mathcal{A}}
\def\B{\mathcal{B}}
\def\C{\mathcal{C}}
\def\I{\mathcal{I}}
\def\N{\mathbb{N}}
\def\P{\mathcal{P}}
\def\Z{\mathbb{Z}}
\def\e{\epsilon}
\def\d{\delta}

\begin{document}

\title{On the Carath\'{e}odory approach to the construction of a measure}

\author{Ivan Werner
\footnote{Email: ivan\_werner@mail.ru}} 
\date{April 10, 2017}
\maketitle

\begin{abstract}\noindent
 The Carath\'{e}odory theorem on the construction of a measure is generalized by replacing the outer measure with an approximation of it and generalizing the Carath\'{e}odory measurability. The new theorem is applied to obtain dynamically defined measures from constructions of outer measure approximations resulting from sequences of measurement pairs consisting of refining $\sigma$-algebras and measures on them which need not be consistent. A particular case when the measurement pairs are given by the action of an invertible map on an initial $\sigma$-algebra and a measure on it is also considered.


\end{abstract}

\tableofcontents

\section{Introduction}

The mathematical endeavor to construct measures, motivated by the need first for the notions of length, area, volume, integral and later for a description of states of stochastic, dynamical and physical systems, has a very long history. It has its brightest point in Lebesgue's groundbreaking work \cite{L} obtaining a countably additive measure on what he called {\it measurable} sets (not exactly what is now usually called {\it Lebesgue measurable} sets). 

Building on the work of Lebesgue, Carath\'{e}odory found an approach to the construction of a countably additive measure \cite{C}, which is very general and convenient for the proofs and applications, by formalizing the notion of the {\it outer measure} and introducing a more restrictive notion of measurability. 
In the case of the Lebesgue outer measure, resulting from a finite, nonnegative and additive set function on an algebra, the class of the {\it Carath\'{e}odory measurable} sets coincides with the Lebesgue one. The Carath\'{e}odory approach, particularly because it does not require an additive set function, found numerous applications (probably, the most prominent one is the construction of a Hausdorff measures \cite{H}) and, in the modern form, is given in every textbook on {\em Measure Theory} (e.g. see Section 1.11 in \cite{B} for a refined presentation of it).

One particular application of the Carath\'{e}odory approach was the construction of equilibrium states for certain random dynamical systems \cite{Wer10} \cite{Wer12}. It was done through a dynamical extension of the Carath\'{e}odory outer measure (in physics and probability, one not always encounters  consistent parts from which a measure describing a state of the system needs to be constructed). However, it turned out that the problem of finding criteria on when such measures are not zero requires further research \cite{Wer13}. All paths taken by the author to obtain lower bounds for such measures and analyze them \cite{Wer15} led to various auxiliary set functions which go beyond outer measures, but have certain three properties which we call an {\it outer measure approximation}. 

In this article, we generalize the Carath\'{e}odory measurability, prove a generalization of the Carath\'{e}odory Theorem for outer measure approximations and develop a general measure theory for such constructions. It naturally extends the classical Measure Theory and can be called {\em Dynamical Measure Theory}. Although our proof of the generalization of the Carath\'{e}odory Theorem is an adaptation of the well-known proof, the obtained result leads to new possibilities for construction of measures. Moreover, the introduced novelty, a "primordial" set function which measures the degree of approximation to an outer measure, opens up a new dimension in the theory, which increases the potential for its further development from its internal logic, which seems to have been in a deadlock since 1914. This article can be seen as a continuation of the work of Carath\'{e}odory \cite{C} and the first part of \cite{Wer15}. 

It is organized as follows. We start with the introduction of the new measure-theoretic language and the proof of the generalization of the Carath\'{e}odory Theorem in Section \ref{ps}. In Section \ref{ddoms}, we construct the dynamically defined outer measures (DDM) in a general setting, from a sequence of {\it measurement pairs}.  The DDMs on the generated $\sigma$-algebra are then obtained,  in Section \ref{ddms}, from the outer measures in the case of a refining, but not necessarily consistent sequence of measurement pairs consisting of $\sigma$-algebras and measures on them. The outer measure approximations are constructed within the same generality in Subsection \ref{ddmomp}. In Section \ref{imc}, we consider the particular case, in which the constructions significantly simplify, when the measurement pairs are generated by an invertible map from an initial $\sigma$-algebra and a measure on it. We provide some examples in Section \ref{es}.

The developed theory is applied in the next article \cite{Wer15} for computation and analysis of various lower bounds for the DDMs in the case when the measurement pairs are generated by an invertible map.

We will use the following notation in this article.  '$f|_{\A}$' will denote the restriction of a function $f$ on a set $\A$, and $\ll$ will denote the absolute continuity relation for set functions. The set of natural numbers $\N$ starts with $1$.

\section{A generalization of the Carath\'{e}odory theorem}\label{ps}

As indicated in the introduction, we will need a generalization of the Carath\'{e}odory Theorem, in order to obtain some measures in this article. We present it in this section, along with the definitions of some new notions which we are going to use.

Let $X$ be a set and $\P(X)$ be the set of all subsets of $X$.

\begin{Definition} 
	We call a collection $\A$ of subsets of $X$  an {\it aggregate} on $X$ iff\\
	$(i)$ $\emptyset\in\A$, and\\
	$(ii)$ $\bigcup_{i\in\N}A_i\in\A$ if $A_i\in\A$ for all $i\in\N$.
\end{Definition} 
\begin{Definition} 
	Let $\A$ be an aggregate. We call $\mu:\A\longrightarrow[0,+\infty]$ an {\it outer measure on }$\A$ iff\\
	$(i)$ $\mu(\emptyset)=0$,\\
	$(ii)$ $\mu(A)\leq\mu(B)$ for all $A,B\in\A$ with $A\subset B$, and\\
	$(iii)$ $\mu$ is countably subadditive, i.e.  
	\[\mu\left(\bigcup\limits_{i\in\N}A_i\right)\leq\sum\limits_{i\in\N}\mu(A_i)\ \ \ \mbox{ for all }(A_i)_{i\in\N}\subset\A.\]   
	If $\A=\P(X)$,  $\mu$ is called an {\it outer measure on} $X$. We call  $(\A,\mu)$ a {\it measurement pair} on $X$ iff  $\A$ is an aggregate on $X$ and $\mu$ is an outer measure on $\A$. We call an outer measure $\mu$ on $\A$ a {\it measure}  iff it is countably additive, i.e.
	\[\mu\left(\bigcup\limits_{i\in\N}A_i\right)=\sum\limits_{i\in\N}\mu(A_i)\ \ \ \mbox{ for all pairwise disjoint }(A_i)_{i\in\N}\subset\A.\]   
\end{Definition} 

Note that, since $\bigcup_{i=1}^\infty A_i=\bigcup_{i=1}^\infty A_i\setminus(A_{i-1}\cup...\cup A_1)$  where $(A_{i-1}\cup...\cup A_1) := \emptyset$ for $i=1$, it is equivalent to require for the definition of an outer measure that the inequality in (iii) holds true only for pairwise disjoint families of sets, because of (ii), if $\A$ is also a ring.

\begin{Definition}
	We call a set function $\mu$ on an aggregate $\A$ {\it finitely additive} iff
	\[\mu(A\cup B)=\mu(A)+\mu(B)\ \ \ \mbox{ for all disjoint }A,B\in\A.\]
\end{Definition}

Note that, as one easily sees, an outer measure is a measure if and only if it is finitely additive.

\begin{Definition}
	Let $\mu$ be an outer measure on $X$. $A\in\P(X)$ is called {\em Carath\'{e}odory $\mu$-measurable} iff
	\[\mu(Q)=\mu(Q\cap A)+\mu(Q\setminus A)\ \ \ \mbox{ for all }Q\in\P(X).\]
	Let $\A_{\mu}$ denote the class of all Carath\'{e}odory $\mu$-measurable subsets of $X$.
\end{Definition}

In order to formulate the generalization of the Carath\'{e}odory Theorem, we need the following definitions.

Let $\A$ be a $\sigma$-algebra on $X$ and $\nu$ be a non-negative set function on $\A$ such that 
\begin{equation}\label{np}
  \nu\left(\bigcup\limits_{i=1}^\infty A_i\right)\leq\sum\limits_{i=1}^\infty\nu(A_i)<\infty\mbox{ for all pairwise disjoint }(A_i)_{i=1}^\infty\subset\A.
\end{equation}

It appears that the following definition has already been considered as a possible generalization of the Carath\'{e}odory measurability, see Exercise 1.12.150 (p. 102) in \cite{B}, but it seems that it has been dismissed as not leading anywhere. 

\begin{Definition}\label{gCm}
	Let $\mu:\P(X)\longrightarrow [0,+\infty]$ such that $\mu(\emptyset)=0$.
	We call $A\in\P(X)$  {\em Carath\'{e}odory $(\A,\mu)$-measurable} iff
	\[\mu(Q)=\mu(Q\cap A)+\mu(Q\setminus A)\ \ \ \mbox{ for all }Q\in\A.\]
	Let $\A_{\A\mu}$ denote the class of all Carath\'{e}odory $(\A,\mu)$-measurable subsets of $X$.
\end{Definition}

Obviously, $\A_{\A\mu} = \A_{\mu}$ if $\A = \P(X)$.

\begin{Definition}
	Let $f:[0,\infty)\longrightarrow[0,\infty)$ be a non-decreasing function  which is continuous at $0$ with $f(0)=0$. For every $t>0$, let $\mu_t:\P(X)\longrightarrow [0,+\infty]$ be such that $\mu_t\geq \mu_s$ (setwise) for all $t\leq s$, and let $\mu:=\lim_{t\to 0}\mu_t$ (setwise).  We call the family $(\mu_t)_{t>0}$ an {\em outer measure $(\A,\nu,f)$-approximation} iff
	\begin{eqnarray*}
		&(i)&\ \ \mu(\emptyset) = 0,\\
		&(ii)&\ \ \mu_{f(\nu(B\setminus A))+\epsilon}(A)\leq \mu(B)\ \ \ \mbox{ for all }\epsilon>0\mbox{ and }A,B\in\A\mbox{ with }A\subset B\mbox{, and}\\
		&(iii)&\ \ \mu\left(\bigcup\limits_{i=1}^\infty A_i\right)\leq \sum\limits_{i=1}^\infty\mu\left(A_i\right) \mbox{ for all pairwise disjoint }(A_i)_{i=1}^\infty\subset\A.
	\end{eqnarray*}
\end{Definition}

Observe that $\mu$ is an outer measure on $X$ if $\A=\P(X)$ and $\nu(A)=0$ for all $A\in\A$. 

If $(\mu_t)_{t>0}$ is simply a non-decreasing family of outer measures on $X$ (as, for example, in the case of the construction of a Hausdorff measure), then $\mu$ is automatically an outer measure on $X$, and Definition \ref{gCm} does in fact seem not to lead anywhere (see Exercise 1.12.150 (p. 102) in \cite{B}). In general however (the examples of such $(\mu_t)_{t>0}$ are given in Subsection \ref{ddmomp}), it leads to the following theorem.

\begin{theo}\label{ctg}
	Suppose $f:[0,\infty)\longrightarrow[0,\infty)$ is a non-decreasing function  which is continuous at $0$ with $f(0)=0$,  $(\mu_t)_{t>0}$ is an outer measure $(\A,\nu,f)$-approximation and $\mu:=\lim_{t\to 0}\mu_t$. Then $\A\cap\A_{\A\mu}$ is a $\sigma$-algebra, and the restriction of $\mu$ on $\A\cap\A_{\A\mu}$ is a measure.
\end{theo}
{\it Proof.}
The proof is an adaptation of the proof of Theorem 5.3 in \cite{Ba}.

Since, by the definition, $A\in\A_{\A\mu}$ iff
\begin{equation}\label{me}
\mu(Q) = \mu(Q\cap A)+\mu(Q\setminus A)\mbox{ for all }Q\in\A,
\end{equation}
we see that $X\in\A_{\A\mu}$ and, by the symmetry, $X\setminus A \in\A_{\A\mu}$ for every $A\in\A_{\A\mu}$. In particular,  $X\setminus A \in\A\cap\A_{\A\mu}$ for every $A\in\A\cap\A_{\A\mu}$.

Let $A,B\in\A\cap\A_{\A\mu}$. We show now that $A\cup B\in\A\cap\A_{\A\mu}$. Replacing $Q$ in \eqref{me} with $Q\cap B$ and $Q\setminus B$ gives two equations the summation of which gives 
\begin{equation}\label{se}
\mu(Q) = \mu(Q\cap A\cap B) + \mu(Q\cap B\setminus A)+\mu(Q\cap A\setminus B) + \mu(Q\setminus (A\cup B))
\end{equation} 
for all $Q\in\A$. Now, replacing $Q$ in \eqref{se} with $Q\cap(A\cup B)$ gives
\begin{equation}\label{te}
\mu(Q\cap(A\cup B)) = \mu(Q\cap A\cap B) + \mu(Q\cap B\setminus A)+\mu(Q\cap A\setminus B) 
\end{equation} 
for all $Q\in\A$. The latter together with \eqref{se} implies that
\[\mu(Q) = \mu(Q\cap (A\cup B)) + \mu(Q\setminus (A\cup B))\ \ \ \mbox{for all }Q\in\A.\]
That is $A\cup B\in\A_{\A\mu}$, and therefore, $A\cup B\in\A\cap\A_{\A\mu}$.

Now, let $(A_i)_{i=1}^\infty\subset \A\cap\A_{\A\mu}$ be pairwise disjoint. Then setting $A=A_1$ and $B=A_2$ in \eqref{te} gives
\[\mu(Q\cap(A_1\cup A_2)) = \mu(Q\cap A_1) + \mu(Q\cap A_2)\ \ \ \mbox{ for all }Q\in\A.\]
Hence, by the induction,
\begin{equation}\label{fme}
\mu\left(Q\cap\bigcup\limits_{i=1}^nA_i\right) = \sum\limits_{i=1}^n\mu\left(Q\cap A_i\right) \ \ \ \mbox{ for all }Q\in\A\mbox{ and }n\geq 1.
\end{equation} 
Let us abbreviate $C_n:=\bigcup_{i=1}^nA_i$, $n\geq 1$, and $C:=\bigcup_{i=1}^\infty A_i$. Then, by the above, $C_n\in \A\cap\A_{\A\mu}$ for all $n\geq 1$, and $C\in \A$. 
Observe that $Q\setminus C\subset Q\setminus C_n$ and $(Q\setminus C_n)\setminus(Q\setminus C) = (Q\cap C)\setminus C_n$ for all $n\geq 1$. Hence, for every $Q\in\A$ and $n\geq 1$, by the property (ii) of the outer measure $(\A,\nu,f)$-approximation, \eqref{fme}, \eqref{np} and the monotonicity of $f$,
\begin{eqnarray*}
	\mu(Q) &=& \mu(Q\cap C_n)+\mu(Q\setminus C_n)\\
	&\geq&\mu(Q\cap C_n)+\mu_{f(\nu((Q\cap C)\setminus C_n))+1/n}(Q\setminus C)\\
	&\geq&\sum\limits_{i=1}^n\mu\left(Q\cap A_i\right)+\mu_{f\left(\sum\limits_{i=n+1}^\infty\nu(Q\cap A_i)\right)+1/n}(Q\setminus C).\\
\end{eqnarray*}
Therefore, by \eqref{np}, since  $f$ is continuous at $0$,
\begin{equation}\label{fe}
\mu(Q)\geq \sum\limits_{i=1}^\infty\mu\left(Q\cap A_i\right)+\mu(Q\setminus C) \ \ \ \mbox{ for all }Q\in\A.
\end{equation}
Hence, by the property (iii) of the outer measure $(\A,\nu,f)$-approximation,
\[\mu(Q)\geq \mu\left(Q\cap \bigcup\limits_{i=1}^\infty A_i\right)+\mu\left(Q\setminus  \bigcup\limits_{i=1}^\infty A_i\right) \ \ \ \mbox{ for all }Q\in\A.\]
Since, by the property (iii) of the outer measure $(\A, \nu,f)$-approximation, holds true also the reverse inequality, $\bigcup_{i=1}^\infty A_i\in\A_{\A\mu}$, and therefore, $\bigcup_{i=1}^\infty A_i\in\A\cap\A_{\A\mu}$. Hence, the algebra $\A\cap\A_{\A\mu}$ is a $\cap$-stable  Dynkin system, and therefore, it is a $\sigma$-algebra.

Finally, putting $Q= \bigcup_{i=1}^\infty A_i$ in \eqref{fe} and the property (iii) of the outer measure $(\nu,\A,f)$-approximation give 
\[\mu\left(\bigcup\limits_{i=1}^\infty A_i\right)= \sum\limits_{i=1}^\infty\mu\left(A_i\right).\]
Thus $\mu$ is a measure on $\A\cap \A_{\A\mu}$.
\hfill$\Box$

Clearly, Theorem \ref{ctg} reduces to the Carath\'{e}odory Theorem if $\A=\P(X)$ and $\nu(A)=0$ for all $A\in\A$. 

The theorem opens the question on the existence of even more general approximating families of set functions and measurability definitions which also lead to measures. In particular, the reader might find it curious that $\nu$ does not play any role in Definition \ref{gCm}. Hopefully, the structure of the theory is now rich enough to drive its further development from its internal logic.

\section{The dynamically defined outer measure}\label{ddoms}

Now, we define, if not a proper generalization, then at least  a dynamical extension of the Carath\'{e}odory outer measure, with a particular case of which this article is concerned.

Since the main problem with the measures obtained in \cite{Wer10} is to determine when they are not zero, naturally arises the question whether one can admit for the construction also a sequence of measurement pairs with increasing norms. 

We will use the opportunity presented in this paper to explore also the question on how far the generalization can be pushed. 

Let $\I$ be a countable set and $(\A_m,\phi_m)_{m\in \I}$ be a family of measurement pairs on $X$.
\begin{Definition} 
	For $Q\in P(X)$, define
	\[\C(Q):=\left\{(A_m)_{m\in\I}|\ A_m\in\A_{m}\ \forall m\in\I\mbox{ and }Q\subset\bigcup\limits_{m\in\I}A_m \right\}\] 
	and
	\begin{equation*}
	\Phi(Q):=\left\{\begin{array}{cc}
	\inf\limits_{(A_m)_{m\in\I}\in\C(Q)}\sum\limits_{m\in\I}\phi_{m}(A_m) & \mbox{if }\C(Q)\neq\emptyset\\
	+\infty  &\mbox{ otherwise.}
	\end{array}\right.
	\end{equation*}
\end{Definition} 

\begin{lemma}\label{oml}
	$\Phi$ is an outer measure on $X$.
\end{lemma}
{\it Proof.}
Clearly, $\Phi(\emptyset)=0$.

Let $Q_1\subset Q_2\subset X$. Then, obviously, $\C(Q_2)\subset\C(Q_1)$. Hence
\[\Phi(Q_1)\leq\Phi(Q_2).\] 

Let $Q_n\subset X$ for all $n\in\N$ and $\epsilon>0$. Clearly, for the proof of the countable subadditivity, we can assume that $\sum_{n\in\N}\Phi(Q_n)<\infty$. Then, for each $n\in\N$, there exists $(A_{nm})_{m\in\I}\in\C(Q_n)$ such that
\[\sum\limits_{m\in\I}\phi_m(A_{nm})<\Phi(Q_n)+\epsilon2^{-n}.\]
Then $(\bigcup_{n\in\N}A_{nm})_{m\in\I}\in\C(\bigcup_{n\in\N}Q_n)$, and therefore,
\begin{eqnarray*}
	\Phi\left(\bigcup\limits_{n\in\N}Q_n\right)&\leq&\sum\limits_{m\in\I}\phi_m\left(\bigcup_{n\in\N}
	A_{nm}\right)\\
	&\leq&\sum\limits_{m\in\I}\sum\limits_{n\in\N}\phi_m\left(A_{nm}\right)\\
	&\leq&\sum\limits_{n\in\N}\Phi(Q_n)+\epsilon.
\end{eqnarray*}
\hfill$\Box$

\section{The dynamically defined measures (DDM)}\label{ddms}

In this section, we introduce some additional conditions on the measurement pairs which allow to obtain useful measures from the dynamically defined outer measure.

Let $(\A_m,\phi_m)_{m\in \Z\setminus\N}$ be a sequence of measurement pairs on $X$  such that $\mathcal{A}_0\subset\mathcal{A}_{-1}\subset\mathcal{A}_{-2} ...$ Let $\B$ denote the $\sigma$-algebra generated by $\bigcup_{m\leq 0}\mathcal{A}_m$.

\begin{Definition}\label{dom}
	For every $Q\in\P(X)$ and $i\in\Z\setminus\N$, define
	\[\C_i(Q):=\left\{(A_m)_{m\leq 0}|\ A_m\in\mathcal{A}_{m+i}\
	\forall m\leq0\mbox{ and }Q\subset\bigcup\limits_{m\leq 0}A_m \right\},\]  
	\[\C(Q):= \C_0(Q),\]
	\begin{equation}\label{iddmd}
	\Phi_{i}(Q):=\left\{\begin{array}{cc}\inf\limits_{(A_m)_{m\leq 0}\in\C_i(Q)}\sum\limits_{m\leq 0}\phi_{m+i}(A_m)& \mbox{if }\C_i(Q)\neq\emptyset\\
	+\infty  &\mbox{ otherwise}
	\end{array}\right.,
	\end{equation}
	\begin{equation}\label{isddmd}
	\Phi_{(i)}(Q):=\left\{\begin{array}{cc}\inf\limits_{(A_m)_{m\leq 0}\in\C(Q)}\sum\limits_{m\leq 0}\phi_{m+i}(A_m)& \mbox{if }\C(Q)\neq\emptyset\\
	+\infty  &\mbox{ otherwise}
	\end{array}\right.\ \ \ \mbox{ and}
	\end{equation}
	\[\Phi(Q):=\Phi_0(Q)=\Phi_{(0)}(Q).\]
	By Lemma \ref{oml}, each of $\Phi_{(i)}$ and $\Phi_{i}$ defines an outer measure on $X$. Observe that
	\begin{equation}\label{sm}
	\Phi_{(i)}(Q)\leq \Phi_{(i-1)}(Q), \mbox{  and}
	\end{equation}
	\begin{equation}\label{m}
	\Phi_{i}(Q)\leq \Phi_{i-1}(Q)
	\end{equation}
	for all $i\leq 0$, since $(...,A_{-1},A_{0},\emptyset)\in\mathcal{C}(Q)$ for all $(A_m)_{m\leq 0}\in\mathcal{C}(Q)$, and\\ $(...,B_{-1},B_{0},\emptyset)\in\mathcal{C}_i(Q)$ for all $(B_m)_{m\leq 0}\in\mathcal{C}_{i-1}(Q)$ and $i\leq 0$. Also, it is obvious that
	\[\Phi_{i}(Q)\leq \Phi_{(i)}(Q),\]
	since $\mathcal{C}(Q)\subset\mathcal{C}_i(Q)$ for all $i\leq 0$ . Therefore, we can define
	\[\bar\Phi(Q):=\lim\limits_{i\to-\infty}\Phi_{i}(Q), \mbox{ and}\]
	\[\Phi^*(Q):=\lim\limits_{i\to-\infty}\Phi_{(i)}(Q),\]
	which are also outer measures on $X$, with 
	\begin{equation}\label{omr}
	\Phi(Q)\leq\bar\Phi(Q)\leq \Phi^*(Q)\ \ \ \mbox { for all }Q\in\P(X).
	\end{equation}
\end{Definition}

The following lemma corrects Lemma 2 in \cite{Wer10}.
\begin{Definition}
	For $Q\in\P(X)$, let  $\dot\C(Q)$ denote the set of all $(A_m)_{m\leq 0}\in\C(Q)$ such that $A_i\cap A_j=\emptyset$ for all $i\neq j\leq 0$, and set
	\begin{equation*}
	\dot\Phi(Q):=\left\{\begin{array}{cc}\inf\limits_{(A_m)_{m\leq 0}\in\dot\C(Q)}\sum\limits_{m\leq 0}\phi_{m}(A_m)& \mbox{if }\dot\C(Q)\neq\emptyset\\
	+\infty  &\mbox{ otherwise.}
	\end{array}\right.
	\end{equation*}
\end{Definition}
\begin{lemma}\label{dcl}
	Suppose each $\A_m$ is also a ring. Then 
	\[\dot\Phi(Q)=\Phi(Q)\ \ \ \mbox{ for all }Q\in\P(X).\]
\end{lemma}
{\it Proof.}
Let $Q\in\P(X)$. Obviously,
\begin{equation*}
\dot\Phi(Q) \geq \Phi(Q).
\end{equation*}
Now, let $(A_m)_{m\leq 0}\in\C(Q)$. Set $B_0:=A_0$ and
\begin{equation*}
B_m:=A_m\setminus\left(A_{m+1}\cup...\cup A_0\right)\mbox{ for all }m\leq -1.
\end{equation*}
Then $(B_m)_{m\leq 0}\in\dot\C(Q)$ and $B_m\subset A_m$ for all $m\leq 0$. Hence,
\begin{equation*}
\dot\Phi(Q)\leq \sum\limits_{m\leq 0}\phi_m(B_m)\leq\sum\limits_{m\leq 0}\phi_m(A_m)
\end{equation*}
Therefore,
\begin{equation*}
\dot\Phi(Q)\leq\Phi(Q).
\end{equation*}
\hfill$\Box$

\begin{theo}\label{ddf}
	Suppose each $\A_m$ is a $\sigma$-algebra  and each $\phi_m$ is also finitely additive. Then\\
	(i)  $\Phi_{i}$ is a measure on $\mathcal{A}_i$ for all $i\leq 0$, and\\
	(ii) $\B\subset\A_{\bar\Phi}$ and $\B\subset \A_{\Phi^*}$. In particular, the restrictions of $\bar\Phi$ and $\Phi^*$ on $\B$ are measures.
\end{theo}
{\it Proof.} The proof is an adaptation of a part of the proof of Theorem 1 in \cite{Wer10}.

(i) Let $i\leq 0$, $A\in\mathcal{A}_i$, $Q\subset X$ and $(A_m)_{m\leq 0}\in\mathcal{C}_i(Q)$. Then  $(A_m\cap A)_{m\leq 0}\in\mathcal{C}_i(Q\cap A)$ and $(A_m\setminus A)_{m\leq 0}\in\mathcal{C}_i(Q\setminus A)$. Therefore,
\begin{eqnarray*}
	\sum\limits_{m\leq 0}\phi_{m-i}(A_m)&=&\sum\limits_{m\leq 0}\phi_{m-i}(A_m\cap A)+\sum\limits_{m\leq 0}\phi_{m-i}(A_m\setminus A)\\
	&\geq&\Phi_i(Q\cap A)+\Phi_i(Q\setminus A).
\end{eqnarray*}
Hence,
\begin{equation}\label{min}
\Phi_i(Q)\geq\Phi_i(Q\cap A)+\Phi_i(Q\setminus A).
\end{equation}
Hence, $\mathcal{A}_i\subset\A_{\Phi_{i}}$. Thus the assertion follows by the Carath\'{e}odory Theorem.

(ii) Let $A\in\bigcup_{m\leq 0}\mathcal{A}_m$ and $Q\subset X$. Then there exists $i_0\leq 0$ such that $A\in\mathcal{A}_{m+i}$ for all $m\leq0$ and $i\leq i_0$. For $i\leq i_0$, let $(A_m)_{m\leq 0}\in\mathcal{C}_i(Q)$. Then, as above, we obtain  inequality (\ref{min}), and the limit of the latter gives
\[\bar\Phi(Q)\geq\bar\Phi(Q\cap A)+\bar\Phi(Q\setminus A).\]
Therefore, $\bigcup_{m\leq 0}\mathcal{A}_m\subset\A_{\bar\Phi}$. Since $\B$ is the smallest $\sigma$-algebra containing $\bigcup_{m\leq 0}\mathcal{A}_m$,  it follows by the Carath\'{e}odory Theorem that $\B\subset\A_{\bar\Phi}$ and $\bar\Phi$ is a measure on it. 

Now, turning to  $\Phi^*$, set
\begin{equation*}
B_m:=\left\{\begin{array}{cc}
A_{m-i}\cap A&  \mbox{if }m \leq i\\
\emptyset& \mbox{ otherwise, }
\end{array}\right.
\end{equation*}
and
\begin{equation*}
C_m:=\left\{\begin{array}{cc}
A_{m-i}\setminus A&  \mbox{if }m\leq i\\
\emptyset& \mbox{ otherwise }
\end{array}\right.
\end{equation*}
for all $m\leq 0$.
Then $(B_m)_{m\leq 0}\in\mathcal{C}(Q\cap A)$ and $(C_m)_{m\leq 0}\in\mathcal{C}(Q\setminus A)$. Therefore,
\begin{eqnarray*}
	\sum\limits_{m\leq 0}\phi_{m+2i}(A_m)&=&\sum\limits_{m\leq 0}\phi_{m+2i}(A_m\cap A)+\sum\limits_{m\leq 0}\phi_{m+2i}(A_m\setminus A)\\
	&=&\sum\limits_{m\leq 0}\phi_{m+i}(B_m)+\sum\limits_{m\leq 0}\phi_{m+i}(C_m)\\
	&\geq&\Phi_{(i)}(Q\cap A)+\Phi_{(i)}(Q\setminus A).
\end{eqnarray*}
Hence
\[\Phi_{(2i)}(Q)\geq\Phi_{(i)}(Q\cap A)+\Phi_{(i)}(Q\setminus A).\]
Taking the limit gives
\[\Phi^*(Q)\geq\Phi^*(Q\cap A)+\Phi^*(Q\setminus A).\]
Therefore, $\bigcup_{m\leq 0}\mathcal{A}_m\subset\A_{\Phi^*}$. Thus, by the Carath\'{e}odory Theorem, $\B\subset\A_{\Phi^*}$ and $\Phi^*$ is a measure on it.
\hfill$\Box$

We will denote the measures obtained in Theorem \ref{ddf} also with $\bar\Phi$ and $\Phi^*$ if no confusion is possible. Of these two measures, we will refer to  $\bar\Phi$ as the {\it dynamically defined measure (DDM)}.

\subsection{The DDMs from outer measure approximations}\label{ddmomp}

Observe that $\Phi_k(A)\leq\phi_m(A)$ for all $A\in\A_i$ and $m\leq k\leq i\leq 0$, since $(...,\emptyset,\emptyset,A,\emptyset,...,\emptyset)\in\C_k(A)$. As a result, $\bar\Phi(A)\leq\liminf_{m\to-\infty}\phi_m(A)$. Hence, since we do not assume the consistency of the measurement pairs, the norm of $\bar\Phi$ can be very small or even zero (e.g. see  Example \ref{ex1} for a zero case). Therefore, to make the theory easier to apply, it would be helpful to have some criteria on when a DDM  is not zero.
  
One way towards them, is by relating the inconsistent sequence of  measurement pairs with a consistent one, the existence of which may be known through a non-constructive and less descriptive argument (such as Krylov-Bogolyubov or some other non-unique fixed point theorem). The latter extends to a measure on the generated $\sigma$-algebra through the standard extension procedure (e.g.  Proposition \ref{ip}) and may provide some information on the DDM through some residual relation to it. 

For example, a natural way of relating for this purpose is by obtaining intermediate measures resulting from an integration of some transformations of the density functions with respect to some mutually absolutely continuous measure (e.g. as in Kullback-Leibler divergence, Hellinger integral, etc.), which can be estimated in a particular case and provide a clear residual relation to the original DDM (e.g. through some convex inequality). 

It turns out that there is a general measure-theoretic technique for the construction of such intermediate measures, which naturally extends the dynamically defined outer measure. It allows us even to obtain some computable estimates on $\bar\Phi$ in \cite{{Wer15}}. We develop this technique in this subsection. It requires the generalization of the Carath\'{e}odory Theorem on outer measure approximations proved in Section \ref{ps} (Theorem \ref{ctg}).

Let $(\A_m,\phi_m)_{m\in\Z\setminus\N}$ and $(\A_m,\psi_m)_{m\in\Z\setminus\N}$ be families of measurement pairs on $X$ such that $\mathcal{A}_0\subset\mathcal{A}_{-1}\subset\mathcal{A}_{-2} ...$ and $\bar\Phi$ is finite. (For example, given measure spaces $(\A_m,\phi_m)_{m\in\Z\setminus\N}$ and a measure $\Lambda$ on $\B$ such that $\Lambda\ll\phi_m$ for all $m$, one can consider $\psi_m(A):=\int_{A}(d\Lambda/d\phi_m)^\alpha d\phi_m$ for all $A\in\A_m$, $m\leq 0$ and a fixed $\alpha\in(0,1)$. For further examples, see \cite{Wer15}.)

\begin{Definition}\label{ddma}
	Let $\e>0$, $i\in\Z\setminus\N$ and $Q\in\P(X)$. Let $\C_{\phi,\e,i}(Q)$ denote the set of all $(A_m)_{m\leq 0}\in\C_i(Q)$ such that 
	\[\bar\Phi(Q)>\sum\limits_{m\leq 0}\phi_{m+i}(A_m)-\e,\]
	and abbreviate $\C_{\phi,\e}(Q):=\C_{\phi,\e,0}(Q)$. Define
	\begin{equation}\label{omar}
	\Psi_{\phi,\e,i}(Q):=\inf\limits_{(A_m)_{m\leq 0}\in\C_{\phi,\e,i}(Q)}\sum\limits_{m\leq 0}\psi_{m+i}(A_m)\ \ \ \mbox{ and }\ \ \ \Psi_{\phi,\e}(Q):= \Psi_{\phi,\e,0}(Q).
	\end{equation}
	Observe that $\C_{\phi_0,\d,i}(Q)\subset\C_{\phi,\e,i}(Q)$ for all $0<\delta\leq\e$. Hence,
	\begin{equation}\label{omam}
	\Psi_{\phi,\e,i}(Q)\leq\Psi_{\phi,\delta,i}(Q)\ \ \ \mbox{ for all }0<\delta\leq\e.
	\end{equation}
	Define
	\begin{equation}\label{omarl}
	\Psi_{\phi,i}(Q):=\lim\limits_{\e\to 0}\Psi_{\phi,\e,i}(Q)\ \ \ \mbox{ for all }Q\in\P(X).
	\end{equation}
\end{Definition}

Crucial for our construction is the following property.

\begin{lemma}\label{dmoma}
	\begin{equation*}
	\Psi_{\phi,\e,i}(Q)\leq\Psi_{\phi,\e,i-1}\left(Q\right)\ \ \ \mbox{ for all }Q\in\P(X),\ \e>0\mbox{ and }i\in\Z\setminus\N.
	\end{equation*}
\end{lemma}
{\it Proof.}
Let $Q\in\P(X)$, $\e>0$, $i\in\Z\setminus\N$ and $(C_m)_{m\leq 0}\in\C_{\phi,\e, i-1}(Q)$. Set $D_0:=\emptyset$ and $D_m:=C_{m+1}$ for all $m\leq-1$.
Then $(D_m)_{m\leq 0}\in\C_{i}(Q)$, and
\[\sum_{m\leq 0}\phi_{m+i}(D_m) =\sum_{m\leq -1}\phi_{m+i}(C_{m+1}) =\sum_{m\leq 0}\phi_{m+i-1}(C_{m}) <\bar\Phi(Q)+\e.\]
Hence, $(D_m)_{m\leq 0}\in\C_{\phi,\e,i}(Q)$. Therefore, 
\[\Psi_{\phi,\e,i}(Q)\leq\sum\limits_{m\leq0}\psi_{m+i}\left(D_m\right) = \sum\limits_{m\leq-1}\psi_{m+i}\left(C_{m+1}\right) =\sum\limits_{m\leq0}\psi_{m+i-1}\left(C_{m}\right).\]
Thus the assertion follows.
\hfill$\Box$

By  Lemma \ref{dmoma}, we can make the following definitions.

\begin{Definition}\label{siomad}
	For $\e>0$ and $Q\in\P(X)$, set
	\[\bar\Psi_{\phi,\e}(Q):=\lim\limits_{i\to -\infty}\Psi_{\phi,\e,i}(Q), \mbox{ and }\]
	\[\bar\Psi_{\phi}(Q):=\lim\limits_{\e\to 0}\bar\Psi_{\phi,\e}(Q).\]
\end{Definition}

One easily checks that
\begin{equation}\label{dme}
\bar\Psi_{\phi}(Q)=\lim\limits_{i\to -\infty}\Psi_{\phi,i}\left(Q\right)\ \ \ \mbox{ for all }Q\in\P(X).
\end{equation}

In the following, we will always use the capitalization rule  to denote the map $(\A_m,\psi_m)_{m\leq0}\longrightarrow\bar\Psi_{\phi}$, e.g. $\bar\Phi_{\phi}$ denotes the set function \eqref{dme} with $(\phi_m)_{m\leq 0}$ in place of $(\psi_m)_{m\leq 0}$ in \eqref{omar}. (One easily checks that $\bar\Phi_{\phi}(Q)=\bar\Phi(Q)$ for all $Q\in\P(X)$.)

\begin{lemma}\label{omap}
	Suppose each $\A_m$ is a $\sigma$-algebra and each $\phi_m$ is also finitely additive. Let $i\in\Z\setminus\N$. Then 
	$( \Psi_{\phi,\e,i})_{\e>0}$ and $( \bar\Psi_{\phi,\e})_{\e>0}$ are  outer measure $(\A_{\bar\Phi},\bar\Phi, id)$-approximations.
\end{lemma}
{\it Proof.}
The assertion  that $( \bar\Psi_{\phi,\e})_{\e>0}$ is an  outer measure $(\A_{\bar\Phi},\bar\Phi, id)$-
approximation follows from that for $( \Psi_{\phi,\e,i})_{\e>0}$ by Lemma \ref{dmoma}.

Let $\e>0$. Since $(...,\emptyset,\emptyset)\in\C_{\phi,\e,i}(\emptyset)$, $\Psi_{\phi,\e,i}(\emptyset) = 0$ for all $\e>0$. Hence, property (i) of the outer measure $(\A_{\bar\Phi},\bar\Phi,id)$-approximation is satisfied for $( \Psi_{\phi,\e,i})_{\e>0}$.

Let  $A, B\in\A_{\bar\Phi}$ such that $A\subset B$ and  $(A_m)_{m\leq 0}\in\C_{\phi,\e,i}(B)$. Then $\bar\Phi(B)>\sum_{m\leq 0}\phi_{m+i}(A_m)-\e$. Hence,  since $\bar\Phi$ is a finite outer measure on  $X$,  $\bar\Phi(A)>\sum_{m\leq 0}\phi_{m+i}(A_m)-\bar\Phi(B\setminus A)-\e$. As a result, $(A_m)_{m\leq 0}\in\C_{\phi,\e+\bar\Phi(B\setminus A),i}(A)$. Hence,
\[\C_{\phi,\e,i}\left(B\right)\subset\C_{\phi,\e+\bar\Phi\left(B\setminus A\right),i}\left(A\right).\] 
Therefore, by \eqref{omam},
\begin{eqnarray*}
	\Psi_{\phi,i}(B)&\geq&\Psi_{\phi,\e,i}\left(B\right)\geq\Psi_{\phi,\e+\bar\Phi\left(B\setminus A\right),i}\left(A\right).
\end{eqnarray*}
This implies the property (ii) of the outer measure $(\A_{\bar\Phi},\bar\Phi,id)$-approximation.

Let $(Q_n)_{n\in\N}\subset\A_{\bar\Phi}$ be pairwise disjoint. Clearly, for the proof of property (iii) of the outer measure $(\A_{\bar\Phi},\bar\Phi,id)$-approximation, we can assume that $\sum_{n\in\N}\Psi_{\phi,i}(Q_n)$ is finite. Then, for each $n\in\N$, we can choose $(A^n_m)_{m\leq 0}\in\C_{\phi,\e2^{-n},i}(Q_n)$ such that
\begin{equation}\label{spi}
\sum\limits_{m\leq 0}\psi_{m+i}(A^n_m)< \Psi_{\phi,\e2^{-n},i}(Q_n)+\e2^{-n}.
\end{equation}
For each $m\leq 0$, set $B_m:=\bigcup_{n\in\N}A^n_m$. Then $B_m\in\A_{m+i}$ for all $m\leq0$, and $\bigcup_{n\in\N}Q_n\subset\bigcup_{m\leq0}B_m$. Furthermore, since, by the Carath\'{e}odory Theorem, $\bar\Phi$ is a measure on $\A_{\bar\Phi}$,
\[\bar\Phi\left(\bigcup\limits_{n\in\N}Q_n\right)=\sum\limits_{n\in\N}\bar\Phi\left(Q_n\right)\geq\sum\limits_{n\in\N}\sum_{m\leq 0}\phi_{m+i}\left(A^n_m\right)-\e>\sum_{m\leq 0}\phi_{m+i}\left(B_m\right)-2\e.\]
Hence, $(B_m)_{m\leq 0}\in\C_{\phi,2\e,i}(\bigcup_{n\in\N}Q_n)$. Therefore, by \eqref{spi},
\begin{eqnarray*}
	\Psi_{\phi,2\e,i}\left(\bigcup_{n\in\N}Q_n\right)&\leq&\sum\limits_{m\leq 0}\psi_{m+i}(B_m)\leq\sum\limits_{n\in\N}\sum\limits_{m\leq 0}\psi_{m+i}(A^n_m)\\
	&\leq&\sum\limits_{n\in\N}\Psi_{\phi,\e2^{-n},i}\left(Q_n\right)+\e\leq\sum\limits_{n\in\N}\Psi_{\phi,i}\left(Q_n\right)+\e.
\end{eqnarray*}
Thus
\[\Psi_{\phi,i}\left(\bigcup_{n\in\N}Q_n\right)\leq\sum\limits_{n\in\N}\Psi_{\phi,i}(Q_n).\]
\hfill$\Box$

\begin{theo}\label{moma}
	Suppose each $\A_m$ is a $\sigma$-algebra and each $\phi_m$ and $\psi_m$ is also finitely additive.  Then 
	$\B\subset\A_{\bar\Phi}\cap\A_{\A_{\bar\Phi}\bar\Psi_{\phi}}$ (in particular, $\bar\Psi_{\phi}$ is a  measure on $\B$).
\end{theo}
{\it Proof.}
Let $\e>0$, $A\in\bigcup_{m\leq 0}\mathcal{A}_m$ and $Q\in \A_{\bar\Phi}$. Then there exists $i_0\in\Z\setminus\N$ such that $A\in\A_i$, $\bar\Phi\left(Q\setminus A\right)-\Phi_i\left(Q\setminus A\right)<\e$ and $\bar\Phi\left(Q\cap A\right)-\Phi_i\left(Q\cap A\right)<\e$ for all $i\leq i_0$. Let  $(A_m)_{m\leq0}\in\C_{\phi,\e,i}(Q)$ for some $i\leq i_0$. Then $(A_m\cap A)_{m\leq0}\in\C_i(Q\cap A)$ and  $(A_m\setminus A)_{m\leq0}\in\C_i(Q\setminus A)$. Furthermore, since $\bar\Phi$ is a finite outer measure on $X$, 
\begin{eqnarray*}
	\bar\Phi\left(Q\cap A\right)&\geq&\bar\Phi\left(Q\right) -\bar\Phi\left(Q\setminus A\right)\\
	&>&\sum\limits_{m\leq 0}\phi_{m+i}\left(A_m\right)-\e-\bar\Phi\left(Q\setminus A\right)\\
	&=&\sum\limits_{m\leq 0}\phi_{m+i}\left(A_m\cap A\right)-\e\\
	&&+\sum\limits_{m\leq 0}\phi_{m+i}\left(A_m\setminus A\right)-\bar\Phi\left(Q\setminus A\right)\\
	&\geq&\sum\limits_{m\leq 0}\phi_{m+i}\left(A_m\cap A\right)-\e+\Phi_i\left(Q\setminus A\right)-\bar\Phi\left(Q\setminus A\right)\\
	&>&\sum\limits_{m\leq 0}\phi_{m+i}\left(A_m\cap A\right)-2\e.
\end{eqnarray*}
Hence, $(A_m\cap A)_{m\leq 0}\in\C_{\phi,2\e,i}(Q\cap A)$. The same way, one sees that $(A_m\setminus A)_{m\leq 0}\in\C_{\phi,2\e,i}(Q\setminus A)$. Therefore,
\begin{eqnarray*}
	\sum\limits_{m\leq 0}\psi_{m+i}\left(A_m\right) &=&\sum\limits_{m\leq 0}\psi_{m+i}\left(A_m\cap A\right)+\sum\limits_{m\leq 0}\psi_{m+i}\left(A_m\setminus A\right)\\
	&\geq&\Psi_{\phi,2\e,i}\left(Q\cap A\right)+\Psi_{\phi,2\e,i}\left(Q\setminus A\right).
\end{eqnarray*}
Hence,
\begin{eqnarray*}
	\Psi_{\phi,\e,i}\left(Q\right)&\geq&\Psi_{\phi,2\e,i}\left(Q\cap A\right)+\Psi_{\phi,2\e,i}\left(Q\setminus A\right).
\end{eqnarray*}
Taking the limit as $i\to-\infty$ implies that
\[\bar\Psi_{\phi,\e}\left(Q\right)\geq\bar\Psi_{\phi,2\e}\left(Q\cap A\right)+\bar\Psi_{\phi,2\e}\left(Q\setminus A\right).\]
Now, taking the limit as $\e\to 0$ gives
\[\bar\Psi_{\phi}\left(Q\right)\geq\bar\Psi_{\phi}\left(Q\cap A\right)+\bar\Psi_{\phi}\left(Q\setminus A\right).\]
Since property (iii) of the outer measure $(\A_{\bar\Phi},\bar\Phi,id)$-approximation gives the inverse inequality, it follows
that $A\in\A_{\A_{\bar\Phi}\bar\Psi_{\phi}}$. Hence, $\bigcup_{m\leq 0}\mathcal{A}_m\subset\A_{\bar\Phi}\cap\A_{\A_{\bar\Phi}\bar\Psi_{\phi}}$. Thus, by  Lemma \ref{omap} and Theorem \ref{ctg}, $\B\subset\A_{\bar\Phi}\cap\A_{\A_{\bar\Phi}\bar\Psi_{\phi}}$, and  $\bar\Psi_{\phi}$  is a measure on $\B$. 
\hfill$\Box$

\subsubsection{An inductive extension of the construction}\label{incs}

It turns out that an inference of the residual relation of a DDM to a consistent measure often requires several intermediate measures constructed successively. However, such constructions always follow the same measure-theoretic pattern which is given through the natural inductive extension of the construction from Subsection \ref{ddmomp}, which does not result in anything beyond outer measure approximations, and the same generalization of the Carath\'{e}odory theorem applies.  It goes as follows.

Suppose, for each $n\in\N$, $(\A_m,\psi_{n,m})_{m\in\Z\setminus\N}$ is a family of measurement pairs on $X$ where each $\A_m$ is a $\sigma$-algebra and each $\psi_{n,m}$ is also finitely additive. (For example, given measure spaces $(\A_m,\phi_m)_{m\in\Z\setminus\N}$ and a measure $\Lambda$ on $\B$ such that $\Lambda\ll\phi_m$ for all $m$, one can consider $\psi_{n,m}(A):=\int_{A}(d\Lambda/d\phi_m)^{\alpha_n} d\phi_m$ for all $A\in\A_m$, $m\leq 0$ and fixed $(\alpha_n)_{n\in\N}\subset[0,1]$. For more examples, see \cite{Wer15}.)

Suppose $\phi_m$'s are finitely additive such that $\bar\Phi(X)<\infty$. Then we can obtain a measure $\bar\Psi_{1}:=\bar\Psi_{1\phi}$  on $\A_{\bar\Phi}\cap\A_{\A_{\bar\Phi}\bar\Psi_{1}}$ as in the previous subsection, with $(\psi_{1,m})_{m\leq 0}$ in place of $(\psi_{m})_{m\leq 0}$.

\begin{Definition}\label{icd}
	Let $Q\in\P(X)$, $\e>0$ and $i\leq 0$. Set $\C_{1,\e,i}(Q):=\C_{\phi,\e,i}(Q)$. Then for $n\geq 2$, provided $\bar\Psi_{k}(Q)<\infty$ for all $k=1,...,n-1$,   we  can define recursively, 
	\[\C_{n,\e,i}(Q):=\left\{(A_m)_{m\leq 0}\in\C_{n-1,\e,i}(Q)|\ \bar\Psi_{n-1}(Q)>\sum\limits_{m\leq 0}\psi_{n-1,m+i}(A_m)-\e\right\},\]
	\[\Psi_{n,\e,i}(Q):=\inf\limits_{(A_m)_{m\leq 0}\in\C_{n,\e,i}(Q)}\sum\limits_{m\leq 0}\psi_{n,m+i}(A_m),\]
	\[\bar\Psi_{n,\e}(Q):=\lim\limits_{i\to-\infty}\Psi_{n,\e,i}(Q)\ \ \ \mbox{ and}\]
	\[\bar\Psi_{n}(Q):=\lim\limits_{\e\to 0}\bar\Psi_{n,\e}(Q),\]
	since, as one easily verifies the same way as in the proof of Lemma \ref{dmoma},  $\Psi_{n,\e,i}(Q)\leq\Psi_{n,\e,i-1}(Q)$ and, obviously, $\Psi_{n,\e,i}(Q)\leq\Psi_{n,\d,i}(Q)$ for all $i\leq0$ and $0<\d\leq\e$. Let us abbreviate $\Psi_{n,\e}(Q):=\Psi_{n,\e,0}(Q)$,
	and set
	\[\Psi_{n}(Q):=\lim\limits_{\e\to 0}\Psi_{n,\e}(Q).\]
\end{Definition}

The following corollary does the inductive step.

\begin{cor}\label{idmc}
	Let $n\in\N$. Suppose, for each $k\in\{1,...,n\}$,  $\bar\Psi_{k}$, which is given by the above recursive construction, is a finite measure on a $\sigma$-algebra $\B_k$  such that $\B\subset\B_n\subset...\subset\B_1\subset\A_{\bar\Phi}$. Let \[\nu_n(Q):=\max\limits_{1\leq k\leq n}\left\{\bar\Psi_{k}(Q)\right\}\vee\bar\Phi(Q)\ \ \ \mbox{ for all }Q\in\B_n.\] Then\\
	(i) $(\Psi_{n+1,\e, i})_{\e>0}$ for all $i\leq 0$ and $(\bar\Psi_{n+1,\e})_{\e>0}$ are outer measure $(\B_n,\nu_n,id)$-approximations, and\\
	(ii) $\B_n\cap\A_{\B_n\bar\Psi_{n+1}}$ is $\sigma$-algebra such that $\B\subset\B_n\cap\A_{\B_n\bar\Psi_{n+1}}$, and $\bar\Psi_{n+1}$ is a measure on $\B_n\cap\A_{\B_n\bar\Psi_{n+1}}$.
\end{cor}
{\it Proof.}
(i) Checking, the same way (only with a slight nuance in the proof of the property (ii) of the outer measure approximation), the corresponding steps as in the proof of Lemma \ref{omap} verifies (i).

(ii) Clearly, by the hypothesis, for every pairwise disjoint $(Q_i)_{i\in\N}\subset\B_n$,
\[\nu_n\left(\bigcup\limits_{i\in\N}Q_i\right)\leq\sum\limits_{i\in\N}\nu_n\left(Q_i\right)\leq\bar\Phi(X)+\sum\limits_{k=1}^n\bar\Psi_k(X)<\infty.\] 
Hence, by (i) and Theorem \ref{ctg}, $\B_n\cap\A_{\B_n\bar\Psi_{n+1}}$ is $\sigma$-algebra, and $\bar\Psi_{n+1}$ is a measure on it.

Next, we show that $\B\subset\B_n\cap\A_{\B_n\bar\Psi_{n+1}}$, as the proof of it has some nuances to that of Theorem \ref{moma}.  Let $\e_n>0$, $A\in\bigcup_{m\leq 0}\mathcal{A}_m$ and $Q\in \B_n$. Successively choose  $\e_{n}>\e_{n-1}>...>\e_0>0$ such that
\[\bar\Psi_{k,3\e_{k-1}}(Q\setminus A)>\bar\Psi_{k}(Q\setminus A)-\e_{k}\mbox{ and }\bar\Psi_{k,3\e_{k-1}}(Q\cap A)>\bar\Psi_{k}(Q\cap A)-\e_{k}\]
for all $k=n,...,1$.
Then there exists $i_0\in\Z\setminus\N$ such that for all $i\leq i_0$, $A\in\A_i$,
\[\Phi_i\left(Q\setminus A\right)>\bar\Phi\left(Q\setminus A\right)-\e_0\mbox{, }\Phi_i\left(Q\cap A\right)>\bar\Phi\left(Q\cap A\right)-\e_0,\]
and
\[\Psi_{k,3\e_{k-1},i}(Q\setminus A)>\bar\Psi_{k,3\e_{k-1}}(Q\setminus A)-\e_{0}\mbox{, }\Psi_{k,3\e_{k-1},i}(Q\cap A)>\bar\Psi_{k,3\e_{k-1}}(Q\cap A)-\e_{0}\]
for all $k=n,...,1$. 

Let $i\leq i_0$ and $(A_m)_{m\leq0}\in\C_{n+1,\e_0,i}(Q)$. Then $(A_m\cap A)_{m\leq0}\in\C_i(Q\cap A)$ and  $(A_m\setminus A)_{m\leq0}\in\C_i(Q\setminus A)$. 

Now, we show by induction that
\begin{equation}\label{is}
(A_m\cap A)_{m\leq 0}\in\C_{k,3\e_{k-1},i}(Q\cap A)\mbox{ and }(A_m\setminus A)_{m\leq 0}\in\C_{k,3\e_{k-1},i}(Q\setminus A)
\end{equation}
for all $k=1,...,n+1$.  As in the proof of Theorem \ref{moma}, one sees that  $(A_m\cap A)_{m\leq 0}\in\C_{1,2\e_0,i}(Q\cap A)$ and $(A_m\setminus A)_{m\leq 0}\in\C_{1,2\e_0,i}(Q\setminus A)$. Thus the induction beginning holds true.  Suppose \eqref{is} is true for some $k\in\{1,...,n\}$. Observe that, since $\C_{n+1,\e_0,i}(Q)\subset ...\subset\C_{2,\e_0,i}(Q)$, by the choice of $i_0$ and $\e_{k}$,
\begin{eqnarray*}
	\bar\Psi_{k}\left(Q\setminus A\right)&=&\bar\Psi_{k}\left(Q\right) -\bar\Psi_{k}\left(Q\cap A\right)\\
	&>&\sum\limits_{m\leq 0}\psi_{k,m+i}\left(A_m\right)-\e_0-\bar\Psi_{k}\left(Q\cap A\right)\\
	&=&\sum\limits_{m\leq 0}\psi_{k,m+i}\left(A_m\setminus A\right)-\e_0\\
	&&+\sum\limits_{m\leq 0}\psi_{k,m+i}\left(A_m\cap A\right)-\bar\Psi_{k}\left(Q\cap A\right)\\
	&\geq&\sum\limits_{m\leq 0}\psi_{k,m+i}\left(A_m\setminus A\right)-\e_0+\Psi_{k,3\e_{k-1},i}\left(Q\cap A\right)-\bar\Psi_{k}\left(Q\cap A\right)\\
	&>&\sum\limits_{m\leq 0}\psi_{k,m+i}\left(A_m\setminus A\right)-\e_0-\e_0-\e_k.
\end{eqnarray*}
Therefore, $(A_m\setminus A)_{m\leq 0}\in\C_{k+1,3\e_k,i}(Q\setminus A)$.  Analogously, one verifies the symmetrical  part of \eqref{is}  for $k+1$. 

Hence,
\begin{eqnarray*}
	\sum\limits_{m\leq 0}\psi_{n+1,m+i}\left(A_m\right) &=&\sum\limits_{m\leq 0}\psi_{n+1,m+i}\left(A_m\cap A\right)+\sum\limits_{m\leq 0}\psi_{n+1,m+i}\left(A_m\setminus A\right)\\
	&\geq&\Psi_{n+1,3\e_n,i}\left(Q\cap A\right)+\Psi_{n+1,3\e_n,i}\left(Q\setminus A\right),
\end{eqnarray*}
which implies that
\begin{eqnarray*}
	\Psi_{n+1,\e_0,i}\left(Q\right)&\geq&\Psi_{n+1,3\e_n,i}\left(Q\cap A\right)+\Psi_{n+1,3\e_n,i}\left(Q\setminus A\right).
\end{eqnarray*}
Now, taking first the limit as $i\to-\infty$ and then also as $\e_n\to 0$ gives
\[\bar\Psi_{n+1}\left(Q\right)\geq\bar\Psi_{n+1}\left(Q\cap A\right)+\bar\Psi_{n+1}\left(Q\setminus A\right).\]
Since property (iii) of the outer measure $(\B_n,\nu_{n},id)$-approximation gives the inverse inequality, it follows
that $A\in\A_{\B_n\bar\Psi_{n+1}}$. Hence, $\bigcup_{m\leq 0}\mathcal{A}_m\subset\B_n\cap\A_{\B_n\bar\Psi_{n+1}}$. Thus $\B\subset\B_n\cap\A_{\B_n\bar\Psi_{n+1}}$ by  Lemma \ref{omap} and Theorem \ref{ctg}. This completes the proof of (ii).
\hfill$\Box$

Very useful for applications is the following lemma.
\begin{Definition}
	For $n\in\N$, $\e>0$, $i\in\Z\setminus\N$ and $Q\in\P(X)$, let  $\dot\C_{n,\e,i}(Q)$ denote the set of all $(A_m)_{m\leq 0}\in\C_{n,\e,i}(Q)$ such that $A_k\cap A_j=\emptyset$ for all $k\neq j\leq 0$, and define
	\[\dot\Psi_{n,\e,i}(Q):=\inf\limits_{(A_m)_{m\leq 0}\in\dot\C_{n,\e,i}(Q)}\sum\limits_{m\leq 0}\psi_{n,m+i}(A_m)\]
	where the fact that $\dot\C_{n,\e,i}(Q)$ is not empty is clarified in the proof of the following lemma.
\end{Definition}

\begin{lemma}\label{dcoma}
	$\dot\Psi_{n,\e,i}(Q)=\Psi_{n,\e,i}(Q)$ for all $Q\in\P(X)$, $\e>0$ and $i\in\Z\setminus\N$.
\end{lemma}
{\it Proof.}
Let  $Q\in\P(X)$, $\e>0$ and $i\in\Z\setminus\N$. Obviously,
\begin{equation*}
\dot\Psi_{n,\e,i}(Q) \geq \Psi_{n,\e,i}(Q).
\end{equation*}
Now, let $(A_m)_{m\leq 0}\in\C_{n,\e,i}(Q)$. Set $B_0:=A_0$ and
\begin{equation*}
B_m:=A_m\setminus\left(A_{m+1}\cup...\cup A_0\right)\mbox{ for all }m\leq -1.
\end{equation*}
Then, since $(A_m)_{m\leq 0}\in\C_{1,\e,i}(Q)$,
\[\bar\Phi(Q)>\sum\limits_{m\leq0}\phi_{m+i}\left(A_m\right)-\e\geq\sum\limits_{m\leq0}\phi_{m+i}\left(B_m\right)-\e,\]
and therefore, $(B_m)_{m\leq 0}\in\dot\C_{1,\e,i}(Q)$. 
The same way, it follows that  $(B_m)_{m\leq 0}\in\dot\C_{k,\e,i}(Q)$ for all $k=2,...,n$. Hence,
\begin{equation*}
\dot\Psi_{n,\e,i}(Q)\leq \sum\limits_{m\leq 0}\psi_{n,m+i}\left(B_m\right)\leq\sum\limits_{m\leq 0}\psi_{n,m+i}\left(A_m\right).
\end{equation*}
Thus
\begin{equation*}
\dot\Psi_{n,\e,i}(Q)\leq\Psi_{n,\e,i}(Q).
\end{equation*}
\hfill$\Box$

\subsubsection{Some signed DDMs}\label{sddms}

It is useful for obtaining and studying lower bounds for DDMs to have the following extension of the inductive construction on some signed measures, in order to admit some transformations of the density functions with negative values.

Let $(\A_m,\phi_{m})_{m\in\Z\setminus\N}$ and $(\A_m,\psi_{k,m})_{m\in\Z\setminus\N}$, $k\in\{1,...,n\}$, be the families of measurement pairs for $n\in\N$ where each $\A_m$ is a $\sigma$-algebra and each $\phi_{m}$ and $\psi_{n,m}$  is finitely additive such that $\bar\Phi(X)<\infty$, and $\bar\Psi_{k}(X)<\infty$ for all $k=1,...,n-1$, as in Subsection \ref{incs}. (Note that $\bar\Psi_{n}(X)$ does not need to be finite.) 

Now, for each $k\in\{1,...,n\}$, let $c_k\in[0,\infty)$, and define
\[\psi'_{k,m} := \psi_{k,m}-c_k\phi_{m}\ \ \ \mbox{ for all }m\leq 0.\]
(For example,  given measure spaces $(\A_m,\phi_m)_{m\in\Z\setminus\N}$ and a measure $\Lambda$ on $\B$ such that $\Lambda\ll\phi_m$ for all $m$, one can consider $\psi_{1,m}(A):=\int_{A}(d\Lambda/d\phi_m)^{\alpha} d\phi_m$ and $\psi_{2,m}(A):=\int_{A}(d\Lambda/d\phi_m)^\alpha\log(d\Lambda/d\phi_m) d\phi_m+1/(\alpha e)\phi_m(A)$ for all $A\in\A_m$, $m\leq 0$ and a fixed $\alpha\in(0,1]$. See \cite{Wer15} for further examples.)

\begin{Definition}\label{iscd}
	Let $Q\in\P(X)$, $\e>0$ and $i\leq 0$. Define $\C'_{1,\e,i}(Q):=\C_{\phi,\e,i}(Q)$, and, for $n\geq 2$, define recursively
	\[\C'_{n,\e,i}(Q):=\left\{(A_m)_{m\leq 0}\in\C'_{n-1,\e,i}(Q)|\ \bar\Psi'_{n-1}(Q)>\sum\limits_{m\leq 0}\psi'_{n-1,m+i}(A_m)-\e\right\},\]
	\[\Psi'_{n,\e,i}(Q):=\inf\limits_{(A_m)_{m\leq 0}\in\C'_{n,\e,i}(Q)}\sum\limits_{m\leq 0}\psi'_{n,m+i}(A_m),\]
	\[\bar\Psi'_{n,\e}(Q):=\lim\limits_{i\to-\infty}\Psi'_{n,\e,i}(Q)\ \ \ \mbox{ and}\]
	\[\bar\Psi'_{n}(Q):=\lim\limits_{\e\to 0}\bar\Psi'_{n,\e}(Q)\]
	(analogously to Definition \ref{icd}), since, as one easily verifies the same way as in the proof of Lemma \ref{dmoma},  $\Psi'_{n,\e,i}(Q)\leq\Psi'_{n,\e,i-1}(Q)$ and, obviously, $\Psi'_{n,\e,i}(Q)\leq\Psi'_{n,\d,i}(Q)$ for all $i\leq0$ and $0<\d\leq\e$.
	
	Let us abbreviate $\Psi'_{n,\e}(Q):=\Psi'_{n,\e,0}(Q)$ and define $\Psi'_{n}(Q):=\lim\limits_{\e\to 0}\Psi'_{n,\e}(Q)$.
	
	Define $\dot\C'_{1,\e,i}(Q):=\dot\C_{\phi,\e,i}(Q)$ and let $\dot\C'_{n,\e,i}(Q)$ be the set of all $(A_m)_{m\leq 0}\in\C'_{n,\e,i}(Q)$ such that all $A_m$'s are pairwise disjoint. It will be clear from the proof of the next lemma that $\dot\C'_{n,\e,i}(Q)$ is not empty.
	
	Define $\dot\Psi'_{n,\e,i}(Q)$ and $\bar{\dot\Psi}'_{n}(Q)$ the same way as $\Psi'_{n,\e,i}(Q)$ and $\bar\Psi'_{n}(Q)$ respectively where the infinitum in the definition of $\dot\Psi'_{n,\e,i}(Q)$ is taken over $\dot\C'_{n,\e,i}(Q)$.
	
	Let $\C_{n,\e,i}(Q)$ and $\bar\Psi_{n}(Q)$ be as in Definition \ref{icd}. 
	
	Define $c'_{0}:=0$ and
	\[c'_{n-1}:=\max\limits_{1\leq j\leq n-1}c_j\ \ \ \mbox{ for all }n\geq 2.\]
\end{Definition}

\begin{lemma}\label{sddml}
	(i) \[\Psi_{1}(Q) - c_1\bar\Phi(Q)\leq\Psi'_{1}(Q) \leq  \Psi_{1}(Q) - c_1\Phi(Q)\ \ \ \mbox{ and}\]
	\[\bar\Psi'_{n}(Q) =  \bar\Psi_{n}(Q) - c_n\bar\Phi(Q)\ \ \ \mbox{ for all }Q\in\P(X).\]
	(ii) For every $\e>0$ and $i\leq 0$,
	\[\Psi_{1}(Q) - c_1\bar\Phi(Q)\leq\dot\Psi'_{1}(Q) \leq  \Psi_{1}(Q) - c_1\Phi(Q),\]
	\[\Psi'_{n,\e,i}(Q)\leq\dot\Psi'_{n,\e,i}(Q),\] \[\dot\Psi'_{n,c'_{n-1}(\bar\Phi(Q)-\Phi_i(Q)+\e)+\e,i}(Q)\leq\Psi'_{n,\e,i}(Q)+c_n(\bar\Phi(Q)-\Phi_i(Q)+\e)\ \ \ \mbox{ and}\]
	\[\bar\Psi'_{n}(Q) = \bar{\dot\Psi}'_{n}(Q)\ \ \ \mbox{ for all }Q\in\P(X).\]
\end{lemma}
{\it Proof.} 
(i) The proof is by induction. Let $Q\in\P(X)$, $\e>0$, $i\leq 0$ and  $(A_m)_{m\leq 0}\in\C'_{1,\e,i}(Q)$.  Then, since $\C'_{1,\e,i}(Q) = \C_{1,\e,i}(Q)$,
\begin{eqnarray*}
	\sum\limits_{m\leq 0}\psi'_{1,m+i}\left(S^mA_m\right) + c_1\Phi_i(Q)&\leq&\sum\limits_{m\leq 0}\psi_{1,m+i}\left(S^mA_m\right)\\
	&\leq&\sum\limits_{m\leq 0}\psi'_{1,m+i}\left(S^mA_m\right) + c_1(\bar\Phi(Q)+\e).
\end{eqnarray*}
Therefore,
\[\Psi'_{1,\e,i}(Q)+c_1\Phi_i(Q)\leq \Psi_{1,\e,i}(Q)\leq \Psi'_{1,\e,i}(Q)+ c_1(\bar\Phi(Q)+\e).\]
Thus (i) is true for n = 1. 

Now, suppose we have shown that $\bar\Psi'_{n-j}(Q) =  \bar\Psi_{n-j}(Q) - c_{n-j}\bar\Phi(Q)$ for all $j\in\{1,...,n-1\}$. 
Let $(B_m)_{m\leq 0}\in\C'_{n,\e,i}(Q)$.
Then, for every $j\in\{1,...,n-1\}$, since $(B_m)_{m\leq 0}\in\C'_{n-j,\e,i}(Q)$,
\begin{eqnarray*}
	&&\bar\Psi_{n-j}(Q) - c_{n-j}\bar\Phi(Q)=\bar\Psi'_{n-j}(Q)>\sum\limits_{m\leq 0}\psi'_{n-j,m+i}(B_m)-\e\\
	&=&\sum\limits_{m\leq 0}\psi_{n-j,m+i}(B_m)-c_{n-j}\sum\limits_{m\leq 0}\phi_{m+i}(B_m)-\e\\
	&\geq&\sum\limits_{m\leq 0}\psi_{n-j,m+i}(B_m)-c_{n-j}\left(\bar\Phi(Q) + \e\right)-\e.
\end{eqnarray*}
Hence, for every $j\in\{1,...,n-1\}$,
\[\bar\Psi_{n-j}(Q)>\sum\limits_{m\leq 0}\psi_{n-j,m+i}(B_m) - \left(c_{n-j}+1\right)\e\geq\sum\limits_{m\leq 0}\psi_{n-j,m+i}(B_m) - \left(c'_{n-1}+1\right)\e.\]
Hence, $(B_m)_{m\leq 0}\in\C_{n,\left(c'_{n-1}+1\right)\e,i}(Q)$. (That is 
\begin{equation}\label{scc}
\C'_{n,\e,i}(Q)\subset\C_{n,\left(c'_{n-1}+1\right)\e,i}(Q).)
\end{equation}
Therefore, 
\begin{eqnarray*}
	\Psi_{n,\left(c'_{n-1}+1\right)\e,i}(Q) &\leq&\sum\limits_{m\leq 0}\psi_{n,m+i}(B_m)\leq\sum\limits_{m\leq 0}\psi'_{n,m+i}(B_m)+c_n\sum\limits_{m\leq 0}\phi_{m+i}(B_m)\\
	&\leq&\sum\limits_{m\leq 0}\psi'_{n,m+i}(B_m)+c_n\left(\bar\Phi(Q) + \e\right).
\end{eqnarray*}
Hence,
\begin{equation}\label{smri}
\Psi_{n,\left(c'_{n-1}+1\right)\e,i}(Q) \leq\Psi'_{n,\e,i}(Q)+c_n\left(\bar\Phi(Q) + \e\right).
\end{equation}
Thus
\begin{equation}\label{rsmi}
\bar\Psi_{n}(Q) \leq\bar\Psi'_{n}(Q)+c_n\bar\Phi(Q).
\end{equation}

Now, let $(C_m)_{m\leq 0}\in\C_{n,\e,i}(Q)$. Then, for every $j\in\{1,...,n-1\}$, since $(C_m)_{m\leq 0}\in\C_{n-j,\e,i}(Q)$,
\begin{eqnarray*}
	&&\bar\Psi'_{n-j}(Q) + c_{n-j}\bar\Phi(Q)=\bar\Psi_{n-j}(Q)>\sum\limits_{m\leq 0}\psi_{n-j,m+i}(C_m)-\e\\
	&=&\sum\limits_{m\leq 0}\psi'_{n-j,m+i}(C_m)+c_{n-j}\sum\limits_{m\leq 0}\phi_{m+i}(C_m)-\e\\
	&\geq&\sum\limits_{m\leq 0}\psi'_{n-j,m+i}(C_m)+c_{n-j}\Phi_i(Q) -\e.
\end{eqnarray*}
Hence, for every $j\in\{1,...,n-1\}$,
\begin{eqnarray*}
	\bar\Psi'_{n-j}(Q)&>&\sum\limits_{m\leq 0}\psi'_{n-j,m+i}(C_m) - c_{n-j}\left(\bar\Phi(Q)-\Phi_i(Q)\right)-\e\\
	&\geq&\sum\limits_{m\leq 0}\psi'_{n-j,m+i}(C_m) - c'_{n-1}\left(\bar\Phi(Q)-\Phi_i(Q)\right)-\e.
\end{eqnarray*}
Hence, $(C_m)_{m\leq 0}\in\C'_{n,c'_{n-1}\left(\bar\Phi(Q)-\Phi_i(Q)\right)+\e,i}(Q)$. (That is  
\begin{equation}\label{sco}
\C_{n,\e,i}(Q)\subset\C'_{n,c'_{n-1}\left(\bar\Phi(Q)-\Phi_i(Q)\right)+\e,i}(Q).)
\end{equation}
Therefore, 
\begin{eqnarray*}
	\sum\limits_{m\leq 0}\psi_{n,m+i}(C_m)&=&\sum\limits_{m\leq 0}\psi'_{n,m+i}(C_m)+c_n\sum\limits_{m\leq 0}\phi_{m+i}(C_m)\\
	&\geq&\Psi'_{n,c'_{n-1}\left(\bar\Phi(Q)-\Phi_i(Q)\right)+\e,i}(Q)+c_n\Phi_i(Q).
\end{eqnarray*}
Hence,
\[\Psi_{n,i,\e}(Q)\geq\Psi'_{n,c'_{n-1}\left(\bar\Phi(Q)-\Phi_i(Q)\right)+\e,i}(Q)+c_n\Phi_i(Q).\]
Since there exists $i_0\leq 0$ such that $c'_{n-1}\left(\bar\Phi(Q)-\Phi_i(Q)\right)<\e$ for all $i\leq i_0$, it follows that
\[\Psi_{n,i,\e}(Q)\geq\Psi'_{n,2\e,i}(Q)+c_n\Phi_i(Q)\ \ \ \mbox{ for all }i\leq i_0.\]
Thus, taking the limit as $i\to-\infty$ and then also as $\e\to 0$ implies that
\[\bar\Psi_{n}(Q)\geq\bar\Psi'_{n}(Q)+c_n\bar\Phi(Q),\]
which together with \eqref{smri} and \eqref{rsmi} proves (i).

(ii) Clearly,
\[\Psi'_{1}(Q) \leq \dot\Psi'_{1}(Q)\mbox{ and }\bar\Psi'_{k}(Q) \leq \bar{\dot\Psi}'_{k}(Q)\ \ \ \mbox{ for all }k\in\{1,...,n\}.\]
Hence, by (i),
\[\Psi_{1}(Q) - c_1\bar\Phi(Q)\leq\dot\Psi'_{1}(Q).\]
Define
\begin{equation*}
\begin{array}{cc}
\xi:&\C(Q)\longrightarrow \C(Q)\\
&(A_m)_{m\leq 0}\longmapsto (B_m)_{m\leq 0}
\end{array}
\end{equation*}
by $B_0:=A_0$ and $B_m:=A_m\setminus(A_{m+1}\cup...\cup A_0)$ for all $m\leq-1$.
Let $(A^1_m)_{m\leq 0}\in\C'_{1,\e,i}(Q)$. Set $(B^1_m)_{m\leq 0}:=\xi((A^1_m)_{m\leq 0})$. Then, since $\C'_{1,\e,i}(Q) = \C_{1,\e,i}(Q)$, as in the proof of Lemma \ref{dcoma}, $(B^1_m)_{m\leq 0}\in\dot\C'_{1,\e,i}(Q)$. Therefore,
\begin{eqnarray*}
	\dot\Psi'_{1,\e,i}(Q)&\leq&\sum\limits_{m\leq 0}\psi'_{1,m+i}(B^1_m)\leq\sum\limits_{m\leq 0}\psi_{1,m+i}(B^1_m)-c_1\sum\limits_{m\leq 0}\phi_{m+i}(B^1_m)\\
	&\leq&\sum\limits_{m\leq 0}\psi_{1,m+i}(A^1_m)-c_1\Phi_i(Q).
\end{eqnarray*}
Hence,
\[\dot\Psi'_{1,\e,i}(Q) \leq \Psi_{1,\e,i}(Q)-c_1\Phi_i(Q),\]
and therefore,
\[\dot\Psi'_{1}(Q) \leq \Psi_{1}(Q)-c_1\Phi(Q)\ \ \ \mbox{ and, by (i),}\]
\[\bar{\dot\Psi}'_{1}(Q) \leq \bar\Psi_{1}(Q)-c_1\bar\Phi(Q)=\bar\Psi'_{1}(Q).\]
This proves (ii) for $n=1$.

Now, let $n\geq 2$, $(A^n_m)_{m\leq 0}\in\C'_{n,\e,i}(Q)$ and $(B^n_m)_{m\leq 0}:=\xi((A^n_m)_{m\leq 0})$. Then, for every $1\leq k\leq n-1$,
\begin{eqnarray*}
	\bar\Psi'_{k}(Q)&\geq&\sum\limits_{m\leq 0}\psi'_{k,m+i}(A^n_m) - \e\\
	&=&\sum\limits_{m\leq 0}\psi_{k,m+i}(A^n_m)-c_{k}\sum\limits_{m\leq 0}\phi_{m+i}(A^n_m)-\e\\
	&\geq&\sum\limits_{m\leq 0}\psi_{k,m+i}(B^n_m)-c_{k}\left(\bar\Phi(Q)+\e\right)-\e\\
	&\geq&\sum\limits_{m\leq 0}\psi'_{k,m+i}(B^n_m)+c_{k}\Phi_i(Q)-c_{k}\left(\bar\Phi(Q)+\e\right)-\e\\
	&\geq&\sum\limits_{m\leq 0}\psi'_{k,m+i}(B^n_m)-c'_{n-1}\left(\bar\Phi(Q)-\Phi_i(Q)+\e\right)-\e.
\end{eqnarray*}
Hence,
\[(B^n_m)_{m\leq 0}\in\dot\C'_{n,c'_{n-1}\left(\bar\Phi(Q)-\Phi_i(Q)+\e\right)+\e,i}(Q).\]
Therefore,
\begin{eqnarray*}
	&&\dot\Psi'_{n,c'_{n-1}\left(\bar\Phi(Q)-\Phi_i(Q)+\e\right)+\e,i}(Q)\\
	&\leq&\sum\limits_{m\leq 0}\psi'_{n,m+i}(B^n_m)\\
	&=&\sum\limits_{m\leq 0}\psi_{n,m+i}(B^n_m)-c_n\sum\limits_{m\leq 0}\phi_{m+i}(B^n_m)\\
	&\leq&\sum\limits_{m\leq 0}\psi_{n,m+i}(A^n_m)-c_n\Phi_i(Q)\\
	&\leq&\sum\limits_{m\leq 0}\psi'_{n,m+i}(A^n_m)+c_n\left(\bar\Phi(Q)+\e\right)-c_n\Phi_i(Q).
\end{eqnarray*}
Thus
\[\dot\Psi'_{n,c'_{n-1}(\bar\Phi(Q)-\Phi_i(Q)+\e)+\e,i}(Q)\leq\Psi'_{n,\e,i}(Q)+c_n(\bar\Phi(Q)-\Phi_i(Q)+\e).\]
 In particular, taking successively limits as $i\to-\infty$ and as $\e\to 0$ implies that
\[\bar{\dot\Psi}'_{n}(Q)\leq\bar\Psi'_{n}(Q).\]
This completes the proof of (ii).
\hfill$\Box$

\subsubsection{The consistent case}

In this subsection, we clarify the situation in the important case, on which the majority of contemporary applications of Measure Theory  is based, that of consistent measurement pairs. 

\begin{Definition}
	We call  a family of measurement pairs $(\A_m,\phi_{m})_{m\in\Z\setminus\N}$ on $X$  {\it consistent} iff 
	\begin{equation}\label{ie}
	\phi_m (A) = \phi_{m-1}\left(A\right)\ \ \ \mbox{ for all }A\in\A_m\mbox{ and }m\leq 0.
	\end{equation}
\end{Definition}

If $(\A_m,\phi_{m})_{m\in\Z\setminus\N}$ is consistent, then for every $A\in\bigcup_{m\leq 0}\A_m$ we can define the set function
\[\phi(A):=\phi_m(A)\mbox{ where }m\leq 0\mbox{ such that }A\in\A_m.\]
One easily sees that, because of \eqref{ie}, $\phi$ is well defined and forms a finitely additive measure on the algebra  $\bigcup_{m\leq 0}\A_m$ if  each $\A_m$ is also an algebra and each $\phi_m$ is, in addition, finitely additive, which allows us to connect our construction with the classical results.

In this case, for every $Q\in\P(X)$, define
\[\phi^*(Q):=\inf\left\{\sum\limits_{n\in\N}\phi\left(A_n\right)|\ A_n\in\bigcup\limits_{m\leq 0}\A_m,\ n\in\N,\mbox{ and }Q\subset\bigcup\limits_{n\in\N}\A_n\right\}.\] 
Obviously $\phi^*$ is the usual outer measure introduced by Lebesgue \cite{L} if each $\phi_m$ is finitely additive.

The following proposition is a correction and a generalization of Proposition 1 in \cite{Wer10}.
\begin{prop}\label{ip}
	Suppose $(\A_m,\phi_{m})_{m\in\Z\setminus\N}$ is a consistent family of measurement pairs on $X$ such that each $\A_m$ is a $\sigma$-algebra and each $\phi_m$ is also finitely additive. Then\\
	(i) $\Phi(Q)=\phi^*(Q)=\bar\Phi(Q)$ for all $Q\in\P(X)$, and\\
	(ii) $\Phi(A_m) = \phi_m(A_m)$ for all $A_m\in\A_m$ and $m\leq 0$, and $\Phi$ is the unique extension of $\phi_m$'s on $\B$.
\end{prop}
{\it Proof.}
(i) Let $Q\in\P(X)$. Let $i\leq 0$ and $(A_m)_{m\leq 0}\in\C_i(Q)$. Then
\[\sum\limits_{m\leq 0}\phi_{m+i}(A_m)=\sum\limits_{m\leq 0}\phi(A_m)\geq\phi^*(Q).\]
Hence
\[\Phi_i(Q)\geq\phi^*(Q).\]
Let $(A_n)_{n\in\N}\subset\bigcup_{m\leq 0}\A_m$ such that $Q\subset\bigcup_{n\in\N}A_n$. We will now define recursively $(B_m)_{m\leq 0}\in\C_i(Q)$. Clearly, there exists the greatest $m_1\leq 0$ such that $A_1\in\A_{m_1+i}$. Set $B_{m_1}:=A_1$ and $B_m:=\emptyset$ for all $m_1<m\leq 0$. Assuming that, for some $n\in\N$, we have defined $B_m$ for all $m_n\leq m\leq 0$, choose the greatest $m_{n+1}<m_n$ such that $A_{n+1}\in\A_{m_{n+1}+i}$, and set $B_{m_{n+1}}:=A_{n+1}$ and $B_m:=\emptyset$ for all $m_{n+1}<m<m_n$. Obviously, the procedure defines $B_m$ for all $m\leq 0$ with  desired  properties. Hence,
\[\Phi_i(Q)\leq\sum\limits_{m\leq 0}\phi_{m+i}(B_m)=\sum\limits_{n\in\N}\phi_{m_n+i}(A_{n})=\sum\limits_{n\in\N}\phi(A_{n}).\]
Therefore,
\[\Phi_i(Q)\leq\phi^*(Q).\]
Thus
\[\Phi_i(Q)=\phi^*(Q).\]
Since $i\leq 0$ was arbitrary, this proves (i).

(ii) The assertion follows from (i) and the well known fact that $\phi^*$ always extends the finitely additive measure on an algebra from which it results, and that the measure resulting from the restriction of the outer measure on the $\sigma$-algebra generated by the algebra is a unique extension.
\hfill$\Box$

Somewhat surprisingly, the same can be proved for $\bar\Psi'_{n}$ from Subsection \ref{sddms} if it is finite and non-negative. It is crucial for some estimations of $\Phi(X)$ in \cite{Wer15}.

\begin{prop}\label{coma}
	For $n\in\N$, let $(\A_m,\phi_{m})_{m\in\Z\setminus\N}$ and $(\A_m,\psi'_{k,m})_{m\in\Z\setminus\N}$, $k\in\{1,...,n\}$, be the families of (signed) measurement pairs from Subsection \ref{sddms} such that $\psi'_{n,m} = \psi_{n,m}$ for all $m\leq 0$ and $(\A_m,\psi_{n,m})_{m\in\Z\setminus\N}$ is consistent such that $\bar\Psi'_{n}(X)<\infty$. Then
	
	(i) $\Psi'_{n}(Q)=\psi_n^*(Q)=\bar\Psi'_{n}(Q)$ for all $Q\in\B$, and\\
	(ii) $\Psi'_{n}(A_m) = \psi_{n,m}(A_m)$ for all $A_m\in\A_m$ and $m\leq 0$, and $\Psi'_{n}$ is the unique extension of $\psi_{n,m}$'s on $\B$.
\end{prop}
{\it Proof.}
(i)  Let $Q\in\B(X)$, $i\in\Z\setminus\N$, $\e>0$ and $(A_m)_{m\leq 0}\in\C'_{n,\e,i}(Q)$. Then
\[\sum\limits_{m\leq 0}\psi_{n,m + i}(A_m)\geq\psi_n^*(Q),\mbox{ and therefore,} \]
\begin{equation}\label{fce}
\bar\Psi'_{n}(Q)\geq\Psi'_{n,i}(Q)\geq\psi_n^*(Q).
\end{equation}
On the other hand, by Proposition \ref{ip}, $\psi_n^*$ is a measure on $\B$, which uniquely extends all $\psi_{n,m}$.  Then, for $(B_m)_{m\leq 0}\in\dot\C'_{n,\e,i}(X)$,
\[\sum\limits_{m\leq 0}\psi_{n,m + i}(B_m)=\psi_n^*(X),\mbox{ and therefore, } \]
\[\dot\Psi'_{n,i}(X)\leq\psi_n^*(X).\]
Hence, by Lemma \ref{sddml} (ii),
\[\bar\Psi'_{n}(X)\leq\psi_n^*(X).\]
This, together with \eqref{fce}, implies that $\bar\Psi'_{n}(X)=\Psi'_{n,i}(X)=\psi_n^*(X)$. Thus, since, by Corollary \ref{idmc} (ii) and Lemma \ref{sddml} (i), $\bar\Psi'_{n}$ is also a measure on $\B$, its finiteness and \eqref{fce} imply that
\[\bar\Psi'_{n}(Q)=\Psi'_{n,i}(Q)=\psi_n^*(Q)\ \ \ \mbox{ for all }Q\in\B\mbox{ and }i\leq 0.\]

(ii) It follows from (i), the same way as in the proof of Proposition \ref{ip} (ii).
\hfill$\Box$

The following proposition shows that a consistent sequence of finite and non-negative measurement pairs can be always put in front of the construction from Subsection \ref{incs} without changing the obtained DDMs if they all are finite.
In particular, it demonstrates that the construction from Subsection \ref{ddmomp} is a generalization of the construction of $\bar{\Phi}$.

\begin{prop}
	Let $n\in\N$. Let $(\A_m,\phi_m)_{m\in\Z\setminus\N}$ and $(\A_m,\psi_{k,m})_{m\in\Z\setminus\N, k\in\{1,...,n\}}$ be the sequences of measurement pairs as in Subsection \ref{incs} such that $\bar\Psi_{k}(X)<\infty$ for all $k\in\{1,...,n\}$. Let $(\A_m,\lambda_m)_{m\in\Z\setminus\N}$ be a consistent sequence of finitely additive and finite measurement pairs. Let $\bar\Phi''$,
	$\bar\Psi_1''$,...,$\bar\Psi_{n+1}''$ be the measures on $\B$ given by Definition \ref{icd} applied to $\phi''_m:=\lambda_m$, $\psi_{1,m}'':=\phi_m$,$\psi_{2,m}'':=\psi_{1,m}$,...,\\ $\psi_{n+1,m}'':=\psi_{n,m}$ for all $m\leq 0$. Then
	\[\bar\Phi_\lambda(Q)=\bar\Phi(Q)\mbox{ and }\bar\Psi_{k+1}''(Q)=\bar\Psi_{k}(Q)\ \ \ \mbox{ for all }Q\in\B\mbox{ and }k\in\{1,...,n\}.\]
\end{prop}
{\it Proof.} Let $Q\in\B$. Clearly,
\[\bar\Phi_\lambda(Q)\geq\bar\Phi(Q).\]
Let $\e>0$ and $i\in\Z\setminus\N$.  By Proposition \ref{coma}, there exists $(A_m)_{m\leq 0}\in\C_{\phi,\e,i}(Q)$ such that
\[\bar{\Lambda}(Q)+\e=\bar\Lambda_{\phi}(Q)+\e>\sum\limits_{m\leq 0}\lambda_{m+i}(\A_m).\]
Hence,  $(A_m)_{m\leq 0}\in\C_{\lambda,\e,i}(Q)$, and therefore,
\[\bar{\Phi}(Q)+\e>\sum\limits_{m\leq 0}\phi_{m+i}(\A_m)\geq\Phi_{\lambda,\e,i}(Q).\]
Thus
\[\bar{\Phi}(Q)\geq\bar\Phi_{\lambda}(Q).\]
This proves the first equality and the rest follows analogously by the induction.
\hfill$\Box$

\section{DDMs for invertible maps}\label{imc}

Now, we consider a special case where the measurement pairs $(\A_m,\phi_m)_{m\in\Z\setminus\N}$ are generated by an invertible dynamical system acting on $X$.

Let $S:X\longrightarrow X$ be an invertible map and $\A$ be a $\sigma$-algebra on $X$. [In the following, we will slightly abuse the notation by denoting the map induced by $S$ acting on classes of subsets of $X$  by the same letter.] For $m\in\Z\setminus\N$, let $\A_m$ be the $\sigma$-algebra generated by $\bigcup_{i=m}^\infty S^{-i}\A$ and $\B$ denote the $\sigma$-algebra generated by $\bigcup_{i=-\infty}^\infty S^{-i}\A$. Then, obviously, $\A_m\subset\A_{m-1}\subset\B$ for all $m\leq0$.  Hence, $\B$ from the previous section is contained in this $\B$. On the other hand, since $\bigcup_{i=-\infty}^\infty S^{-i}\A\subset\bigcup_{m\leq 0}\A_m$, one sees that this $\B$ is exactly $\B$ from the previous section.

Furthermore, by considering the class of all $B\in\B$ such that $S^{-1}B\in\B$ and observing that it is a $\sigma$-algebra containing $\bigcup_{i=-\infty}^\infty S^{-i}\A$, one sees that $S$ is $\B$-$\B$-measurable, and, analogously, that the same is true for $S^{-1}$. The same argument with $\A_0$ instead of $\B$ shows that $S$ is also $\A_0$-$\A_0$-measurable.

Let $m\leq 0$, then, since $S^{-m}\bigcup_{i=0}^\infty S^{-i}\A\subset S^{-m}\A_0$, $\A_m\subset S^{-m}\A_0$. On the other hand, by considering the class of all $A\in\A_0$ such that $S^{-m}A\in\A_m$ and observing that it is a $\sigma$-algebra containing $\bigcup_{i=0}^\infty S^{-i}\A$, one sees that $S^{-m}\A_0\subset\A_m$. Hence, 

\begin{equation}\label{ads}
\A_m = S^{-m}\A_0\ \ \ \mbox{ for all }m\leq 0.
\end{equation}

Now, let $\phi_0$ be an  outer measure on $\mathcal{A}_0$.  Define  
\[\phi_{m}:=\phi_{0}\circ S^{m}\ \ \ \mbox{ for all }m\leq0.\] 
Then, clearly, $(\mathcal{A}_m, \phi_m)$ is a measurement pair for every $m\leq0$.
Observe that, for every $i\leq 0$ and $Q\in\P(X)$, $(S^{i}A_m)_{m\leq 0}\in\mathcal{C}(S^{i}Q)$ if $(A_m)_{m\leq 0}\in\mathcal{C}_i(Q)$, and $(S^{-i}A_m)_{m\leq 0}\in\mathcal{C}_i(Q)$ if $(A_m)_{m\leq 0}\in\mathcal{C}(S^{i}Q)$. This implies that
\begin{equation}\label{dd}
\Phi_i (Q) = \Phi\left(S^{i}Q\right)
\end{equation}
for all $i\leq 0$ and $Q\subset X$. Therefore, the outer measure $\bar\Phi$ is $S$-invariant.  

Let $n\in\N$ and $(\psi_k)_{k=1}^n$ be an additional family of outer measures on $\mathcal{A}_0$ such that $\psi_k(X)<\infty$ for all $k\in\{1,...,n-1\}$. For $Q\in\P(X)$, $i\in\Z\setminus\N$ and $\e>0$, let $\C_{n,\e,i}(Q)$ and $\Psi_{n,\e,i}(Q)$ be defined as in Definition \ref{icd} resulting from $(\psi_k\circ S^m)_{m\leq 0}$ for all $k\in\{1,...,n\}$.  Further, we will use the abbreviations $\C_{n,\e}(Q):=\C_{n,\e,0}(Q)$ and $\Psi_{n,\e}(Q):= \Psi_{n,\e, 0}(Q)$.

\begin{lemma}\label{omasl} Let $Q\in\P(X)$, $i\leq 0$ and $\e>0$. Then
	\begin{equation*}
	\Psi_{n,\e,i}(Q)=\Psi_{n,\e}\left(S^iQ\right).
	\end{equation*}
\end{lemma}
{\it Proof.} The proof is by induction. Let $(A_m)_{m\leq 0}\in\C_{1,\e,i}(Q)$. By the $S$-invariance of $\bar\Phi$, one easily sees that $(S^{i}A_m)_{m\leq 0}\in\C_{1,\e}(S^{i}Q)$. This implies that
\[\Psi_{1,\e}\left(S^iQ\right)\leq\Psi_{1,\e,i}(Q).\]
Then observing that $(S^{-i}B_m)_{m\leq 0}\in\C_{1,\e,i}(Q)$ if $(B_m)_{m\leq 0}\in\C_{1,\e}(S^{i}Q)$ implies that
\[\Psi_{1,\e}\left(S^iQ\right)\geq\Psi_{1,\e,i}(Q).\]
This proves the assertion for $n=1$.

Now, suppose we have shown that $\Psi_{k,\e,i}(Q)=\Psi_{k,\e}\left(S^iQ\right)$ for all $k\in\{1,...,n-1\}$.  Then $\bar\Psi_k$ is $S$-invariant for all $k\in\{1,...,n-1\}$. Let  $(C_m)_{m\leq 0}\in\C_{n,\e,i}(Q)$. Then
\[\bar\Psi_k(S^iQ)=\bar\Psi_k(Q)>\sum\limits_{m\leq 0}\psi_{k,m+i}(C_m)-\e=\sum\limits_{m\leq 0}\psi_{k,m}(S^iC_m)-\e\]
for all $k\in\{1,...,n-1\}$. Hence, $(S^{i}C_m)_{m\leq 0}\in\C_{n,\e}(S^{i}Q)$. Therefore,
\[\Psi_{n,\e}(S^iQ)\leq\sum\limits_{m\leq 0}\psi_{n,m}\left(S^{i}C_m\right)=\sum\limits_{m\leq 0}\psi_{n,m+i}\left(C_m\right).\]
Hence,
\[\Psi_{n,\e}(S^iQ)\leq\Psi_{n,\e,i}(Q).\]
Let $(D_m)_{m\leq 0}\in\C_{n,\e}(S^iQ)$. Then, the same way, one sees that $(S^{-i}D_m)_{m\leq 0}\in\C_{n,\e,i}(Q)$, which implies that
\[\Psi_{n,\e}(S^iQ)\geq\Psi_{n,\e,i}(Q).\]
This completes the proof.
\hfill$\Box$

Since, in this case, the sequence $(\phi_m)_{m\leq 0}$ is completely determined by $\phi_0$, we will use the notation $\C_{\phi_0,\e}(Q):=\C_{\phi,\e}(Q)$ and $\Psi_{\phi_0,\e}(Q):=\Psi_{\phi,\e}(Q)$, to indicate that.

It turns out, as the next theorem shows, that the construction of the DDMs can be simplified  in this case.

\begin{lemma}\label{sim}
	$S$ is $\A_{\bar\Phi}$-$\A_{\bar\Phi}$-measurable.
\end{lemma}
{\it Proof.}
Let $A\in\A_{\bar\Phi}$ and $Q\in\P(X)$. Since $\bar\Phi$ is $S$-invariant, 
\[\bar\Phi(Q)=\bar\Phi\left(SQ\right)=\bar\Phi\left(SQ\cap A\right)+\bar\Phi\left(SQ\setminus A\right)=\bar\Phi\left(Q\cap S^{-1}A\right)+\bar\Phi\left(Q\setminus S^{-1}A\right).\]
Thus $ S^{-1}A\in\A_{\bar\Phi}$.
\hfill$\Box$

The following theorem is a generalization of Theorem 1 in \cite{Wer10} (and the proof of it is an adaptation of a part of the proof of the latter).
\begin{theo}\label{ddmt}
	(i) Suppose $\phi_0$ is finitely additive such that $\Phi(X)<\infty$. Then $\Phi(B)=\bar\Phi(B)$ for all $B\in\A_{\bar\Phi}$. In particular,  $\Phi$ and $\Phi^*$ are $S$-invariant measures on $\B$.\\
	(ii) Suppose $\phi_0$ and $\psi_{1,0},...,\psi_{n,0}$ are finitely additive such that $\Phi(X)<\infty$ and, for each $k=1,...,n$, $\bar\Psi_{k}$, which is given by Definition \ref{icd}, is a finite measure on the $\sigma$-algebra $B_k$ given by $\B_1:=\A_{\bar\Phi}\cap\A_{\A_{\bar\Phi}\bar\Psi_{1}}$ and
	$\B_k := \B_{k-1}\cap\A_{\B_{k-1}\bar\Psi_{k-1}}$ for all $k>1$ by Corollary \ref{idmc} (ii). Then $\Psi_{k}(Q)=\bar\Psi_{k}(Q)$ for all $Q\in\B_{k}$ and $k=1,...,n$. In particular, in this case, each $\Psi_{k}$ is a $S$-invariant measure on $\B$.
\end{theo}
{\it Proof.}
(i) Let $B\in\A_{\bar\Phi}$. Then, by \eqref{omr},
\[\Phi(B)\leq\bar\Phi(B).\]
Since, by \eqref{dd}, the restriction of $\bar\Phi$ on $\A_{\bar\Phi}$ is a measure such that $\bar\Phi(X) = \Phi(X)$ and $\Phi$ is an outer measure on $X$,
\begin{eqnarray*}
	\bar\Phi(X\setminus B) = \bar\Phi(X) - \bar\Phi(B)
	\leq\Phi(X) - \Phi(B)
	\leq\Phi(X\setminus B).
\end{eqnarray*}
Hence, using  $X\setminus B$ instead of $B$ in the above gives
\[\Phi(B) \geq \bar\Phi(B),\]
which  implies the desired equality. Thus $\Phi$ is a $S$-invariant measures on $\A_{\bar\Phi}$ and, by Theorem \ref{ddf} (ii), on $\B$.

We show now that $\Phi^*$ is $S$-invariant. Let $i\leq 0$. Since $\Phi_{(i)}$ is the outer measure $\Phi$ where the initial measure on $\mathcal{A}_0$ is $\phi_i$ instead of $\phi_0$, by (\ref{m}) and (\ref{dd}),
\[\Phi_{(i)}(Q)\leq\Phi_{(i)}\left(S^{-1}Q\right)\]
for all $Q\subset X$. On the other hand, since $(S^{-1}A_m)_{m\leq 0}\in\mathcal{C}(S^{-1}Q)$ for all $(A_m)_{m\leq 0}\in\mathcal{C}(Q)$,
\[\Phi_{(i)}(S^{-1}Q)\leq\Phi_{(i-1)}(Q)\]
for all $i\leq 0$ and $Q\subset X$. Combining these inequalities and taking the limit gives
\[\Phi^*(S^{-1}Q)=\Phi^*(Q)\]
for all $Q\subset X$. Thus, by Theorem \ref{ddf} (ii),  $\Phi^*$ is a $S$-invariant measure on $\B$.

(ii)   Let $k\in\{1,...,n\}$ and $Q\in\B_k$. Clearly,
\begin{equation}\label{fise}
\Psi_{k}(Q)\leq\bar\Psi_{k}(Q).
\end{equation}
On the other hand, since $\bar\Psi_{k}$ is a measure on $\B_k$ and, by Lemma \ref{omasl}, $\Psi_{k}(X)=\bar\Psi_{k}(X)$, applying \eqref{fise} to $Q':=X\setminus Q$ and using the property (iii) of outer measure approximation $( \Psi_{k,\e})_{\e>0}$, by Corollary \ref{idmc} (i), implies that
\[\bar\Psi_{k}(Q) = \bar\Psi_{k}\left(X\setminus Q'\right)=\bar\Psi_{k}\left(X\right)-\bar\Psi_{k}\left(Q'\right)\leq\Psi_{k}(X) - \Psi_{k}(Q')\leq\Psi_{k}(Q)\]
since  $\Psi_{k}(X)<\infty$.
Hence,
\[\Psi_{k}(Q)=\bar\Psi_{k}(Q).\]
This proves (ii). 
\hfill$\Box$

\subsection{The DDMs on topological spaces}
In this subsection, we show that the definitions of $\Phi$ and $\Psi'_{n}$ are constructive on compact sets in non-pathological cases. This fact is useful for obtaining criteria for the positivity of $\Phi$, see Remark 2 in \cite{Wer15}.

Let $X$ be a Hausdorff topological space. Suppose $S$ is a homeomorphism of $X$ such that the Borel $\sigma$-algebra $\B(X)\subset\B$. 
Let $n\in\N$. Let $(\A_m,\phi_{0}\circ S^m)_{m\in\Z\setminus\N}$ and $(\A_m,\psi'_{k,0}\circ S^m)_{m\in\Z\setminus\N}$, $k\in\{1,...,n\}$, be the families of (signed) measurement pairs as in Subsection \ref{sddms}.
\begin{Definition}
	Let $Q\in\P(X)$ and $\e>0$. Let $\hat\C(Q)$ be the set of all $(A_m)_{m\leq 0}\in\C(Q)$ such that each $A_m$ is open in $X$ and at most finitely many of them are not empty and $\hat\C'_{k,\e}(Q)$ be the set of all $(A_m)_{m\leq 0}\in\C'_{k,\e}(Q)$ such that each $A_m$ is open in $X$ and at most finitely many of them are not empty for all $k\in\{1,...,n\}$. Define
	\begin{equation*}
	\hat\Phi(Q):=\inf\limits_{(A_m)_{m\leq 0}\in\hat\C(Q)}\sum\limits_{m\leq 0}\phi_0\left(S^mA_m\right),
	\end{equation*}
	\begin{equation*}
	\hat\Psi'_{k,\e}(Q):=\inf\limits_{(A_m)_{m\leq 0}\in\hat\C'_{k,\e}(Q)}\sum\limits_{m\leq 0}\psi'_{k,0}\left(S^mA_m\right)\ \ \ \mbox{ and}
	\end{equation*}
	\[\hat\Psi'_{k}(Q):=\lim\limits_{\e\to 0} \hat\Psi'_{k,\e}(Q)\]
	for all $k\in\{1,...,n\}$ for which $\C'_{k,\e}(Q)$ is not empty and set $\hat\Psi'_{k}(Q):=+\infty$ otherwise.
\end{Definition}

\begin{Definition}
	We call a measurement pair $(A_0,\phi_0)$ {\em regular from above} iff $\A_0$ is a ring and for every $A\in\A_0$ and $\e>0$ there exists $O\in\A_0$ such that $O$ is open in $X$, $A\subset O$ and $\phi_0(O\setminus A)<\e$. 
\end{Definition}

\begin{lemma}		
	Suppose the measurement pairs $(\A_0,\phi_0)$ and $(\A_0, \psi_{k,0})_{k\in\{1,...,n\}}$ are regular from above. Let $Q\subset X$ be compact. Then 
	\begin{equation*}
	\Phi(Q)=\hat\Phi\left(Q\right)\ \ \ \mbox{ and }\ \ \ \Psi'_{k}(Q)=\hat\Psi'_{k}\left(Q\right)\mbox{ for all }k\in\{1,...,n\}.
	\end{equation*}
\end{lemma}
{\it Proof.}
First, we show that 
\begin{equation}\label{scamc}
\Phi(Q)=\hat\Phi\left(Q\right)\ \ \ \mbox{ and }\ \ \ \Psi_{k}(Q) = \hat\Psi_{k}\left(Q\right)\mbox{ for all }k\in\{1,...,n\}.
\end{equation}
where the definition of $\hat\Psi_{k}\left(Q\right)$ is given by the particular case of the definition of $\hat\Psi_{k}'\left(Q\right)$ when $c_j=0$ for all $1\leq j\leq k$. Clearly, 
\begin{equation}\label{scamca}
\Phi(Q)\leq\hat\Phi\left(Q\right)\ \ \ \mbox{ and }\ \ \ \Psi_{k}(Q)\leq \hat\Psi_{k}\left(Q\right)\mbox{ for all }k\in\{1,...,n\}.
\end{equation}
Now, let $\e>0$ and $(A_m)_{m\leq 0}\in\C_{n,\e}(Q)$. Since $S$ is a homeomorphism, by the hypothesis on $\phi_0$ and $(\psi_{k,0})_{k\in\{1,...,n\}}$, for every $m\leq 0$, there exists an open $O_m\in\A_m$ such that $A_m\subset O_m$ and
\[\phi_0\left(S^m\left(O_m\setminus A_m\right)\right)<\e2^{-|m|-1}\ \ \ \mbox{ and }\ \ \ \psi_{k,0}\left(S^m\left(O_m\setminus A_m\right)\right)<\e2^{-|m|-1}.\]
Since $Q$ is compact, there exists $m_0\leq 0$ such that $Q\subset\bigcup_{m_0\leq m\leq 0}O_m$. Set $O'_m:=O_m$ for all $m_0\leq m\leq 0$ and $O'_m:=\emptyset$ for all $m<m_0$. Then $(O'_m)_{m\leq 0}\in\hat\C(Q)$, and therefore, by Theorem \ref{ddmt} (i), since $Q\in\B$ (because it is closed),
\begin{eqnarray}\label{ftc}
\hat\Phi(Q)&\leq& \sum\limits_{m\leq 0}\phi_0\left(S^mO'_m\right)\nonumber\\
&\leq&\sum\limits_{m\leq 0}\phi_0\left(S^mO_m\right)\nonumber\\
&\leq& \sum\limits_{m\leq 0}\phi_0\left(S^mA_m\right)+ \sum\limits_{m\leq 0}\phi_0\left(S^m\left(O_m\setminus A_m\right)\right)\nonumber\\
&<&\Phi(Q) +\e+\e.
\end{eqnarray}
Hence, since $\e$ was arbitrary, it follows the first equality of the assertion. Furthermore, by \eqref{ftc}, $(O'_m)_{m\leq 0}\in\hat\C_{1,2\e}(Q)$, and therefore,
\begin{eqnarray*}
	\hat\Psi_{1,2\e}\left(Q\right) &\leq&\sum\limits_{m\leq 0}\psi_{1,0}\left(S^mO'_m\right)\leq\sum\limits_{m\leq 0}\psi_{1,0}\left(S^mO_m\right)\\
	&\leq&\sum\limits_{m\leq 0}\psi_{1,0}\left(S^mA_m\right)+\sum\limits_{m\leq 0}\psi_{1,0}\left(S^m\left(O_m\setminus A_m\right)\right)\\
	&\leq&\Psi_{1}\left(Q\right)+\e+\e.
\end{eqnarray*}
Hence,
\[ \hat\Psi_{1}\left(Q\right)\leq\Psi_{1}(Q)\mbox{ and }(O'_m)_{m\leq 0}\in\hat\C_{2,2\e}(Q).\]
By the induction, one sees, since $(A_m)_{m\leq 0}\in\C_{k,\e}(Q)$ for all $k\in\{1,...,n\}$, that
\[\hat\Psi_{k}\left(Q\right)\leq\Psi_{k}(Q)\mbox{ and }(O'_m)_{m\leq 0}\in\hat\C_{k+1,2\e}(Q)\]
for all $k\in\{1,...,n-1\}$. Therefore, since $(O'_m)_{m\leq 0}\in\hat\C_{n,2\e}(Q)$,
\begin{eqnarray*}
	\hat\Psi_{n,2\e}\left(Q\right) &\leq&\sum\limits_{m\leq 0}\psi_{n,0}\left(S^mO'_m\right)\leq\sum\limits_{m\leq 0}\psi_{n,0}\left(S^mO_m\right)\\
	&\leq&\sum\limits_{m\leq 0}\psi_{n,0}\left(S^mA_m\right)+\sum\limits_{m\leq 0}\psi_{n,0}\left(S^m\left(O_m\setminus A_m\right)\right)\\
	&\leq&\sum\limits_{m\leq 0}\psi_{n,0}\left(S^mA_m\right)+\e.
\end{eqnarray*}
Thus taking the infimum over $\C_{n,\e}(Q)$ and letting $\e\to 0$ implies that
\[\hat\Psi_{n}\left(Q\right)\leq\Psi_{n}(Q).\]
Together with \eqref{scamca}, this completes the proof of \eqref{scamc}.

Now, we show by induction that 
\begin{equation}\label{shm}
  \hat\Psi'_{k}(Q)=\hat\Psi_{k}(Q)-c_k\hat\Phi(Q)\ \ \ \mbox{ for all }k\in\{1,...,n\},
\end{equation}
which together with \eqref{scamc} will imply the assertion. 

Recall that $\hat\C_{k,\e}(Q)$ is defined as the particular case of $\hat\C'_{k,\e}(Q)$ when $c_j=0$ for all $1\leq j\leq k-1$, $2\leq k\leq n$. Let $(B_m)_{m\leq 0}\in\hat\C'_{1,\e}(Q)$ (recall that $\hat\C'_{1,\e}(Q):=\hat\C_{1,\e}(Q)$, and the latter is not empty by the above). Then, since $\hat\Phi(Q)=\Phi(Q)$,
\begin{eqnarray*}
 \hat\Psi_{1,\e}(Q)&\leq&\sum\limits_{m\leq 0}\psi_{1,0}\left(S^mB_m\right)=\sum\limits_{m\leq 0}\psi'_{1,0}\left(S^mB_m\right)+c_1\sum\limits_{m\leq 0}\phi_{0}\left(S^mB_m\right)\\
 &<&\sum\limits_{m\leq 0}\psi'_{1,0}\left(S^mB_m\right)+c_1(\hat\Phi(Q)+\e).
\end{eqnarray*}
Hence,
\[\hat\Psi_{1}(Q)\leq\hat\Psi'_{1}(Q)+c_k\hat\Phi(Q).\]
On the other hand,
\[\sum\limits_{m\leq 0}\psi_{1,0}\left(S^mB_m\right)=\sum\limits_{m\leq 0}\psi'_{1,0}\left(S^mB_m\right)+c_1\sum\limits_{m\leq 0}\phi_{0}\left(S^mB_m\right)\geq\hat\Psi'_{1,\e}(Q)+c_k\hat\Phi(Q),\]
which implies \eqref{shm} for the case $k=1$. 

Suppose we have proved \eqref{shm} for $k\in\{1,...,n-1\}$. Then, by Lemma \ref{sddml} (i), \eqref{scamc} and the hypothesis, $\hat\C'_{n,\e}(Q)$ is not empty. Let $(C_m)_{m\leq 0}\in\hat\C'_{n,\e}(Q)$. By \eqref{scc}, $(C_m)_{m\leq 0}\in\hat\C_{n,\left(c'_{n-1}+1\right)\e}(Q)$. Therefore,
\begin{eqnarray*}
	\hat\Psi_{n,\left(c'_{n-1}+1\right)\e}(Q)&\leq&\sum\limits_{m\leq 0}\psi_{n,0}\left(S^mC_m\right)\\
	&=&\sum\limits_{m\leq 0}\psi'_{n,0}\left(S^mC_m\right)+c_n\sum\limits_{m\leq 0}\phi_{0}\left(S^mC_m\right)\\
	&\leq&\sum\limits_{m\leq 0}\psi'_{n,0}\left(S^mC_m\right)+c_n\left(\hat\Phi(Q)+\e\right).
\end{eqnarray*}
This implies
\[\hat\Psi_{n}(Q)\leq\hat\Psi'_{n}(Q)+c_n\hat\Phi(Q).\]
Now, let $(D_m)_{m\leq 0}\in\hat\C_{n,\e}(Q)$, where the latter is not empty by the above. By \eqref{sco}, $(D_m)_{m\leq 0}\in\hat\C'_{n,\e}(Q)$. Therefore,
\[\sum\limits_{m\leq 0}\psi_{n,0}\left(S^mD_m\right)=\sum\limits_{m\leq 0}\psi'_{n,0}\left(S^mD_m\right)+c_n\sum\limits_{m\leq 0}\phi_{0}\left(S^mD_m\right)\geq\hat\Psi'_{n,\e}(Q)+c_n\hat\Phi(Q),\]
which implies the converse inequality and completes the proof of \eqref{shm}.
\hfill$\Box$

\subsection{The norm of the DDM and the non-invariance of the initial measure}

The next proposition  states clearly the obvious dependence of the norm of the DDM on how far the  initial measure is from being invariant. It is a generalization of Proposition 2 in \cite{Wer10}.

\begin{prop}\label{nnp}
	Suppose $\phi_0$ is finitely additive such that $\phi_0(X)<\infty$. Then\\
	(i) \[\Phi(X)\leq\phi_0(X) -\sup\limits_{m\leq 0}\sup\limits_{A\in\A_0}\left|\phi_0\left(S^mA\right)-\phi_0(A)\right|,\mbox{ and}\]
	(ii)  the following are equivalent:\\
	a) $\Phi(X)=\phi_0(X)$,\\
	b) $\phi_0(S^{-1}A)=\phi_0(A)$ for all $A\in\A_0$.
\end{prop}
{\it Proof.}
Let $k,m\leq 0$ and $A\in\A_0$. Then, since, by Theorem \ref{ddmt} (i), $\Phi$ is a $S$-invariant measure on $\B$,
\[\phi_0(X)-\phi_0\left(S^mA\right)=\phi_0\left(X\setminus S^mA\right)\geq\Phi\left(X\setminus A\right)=\Phi(X)-\Phi(S^kA).\] 
Hence,
\[\phi_0(X)-\Phi(X)\geq\phi_0\left(S^mA\right)-\Phi(S^kA)\geq\phi_0\left(S^mA\right)-\phi_0(S^kA).\]
Thus (i) follows.

(ii) The implication from a) to b) follows by (i). The converse follows by Proposition \ref{ip} (ii).
\hfill$\Box$

\subsection{The absolute continuity of the DDMs}

The following lemma is the first piece which can be salvaged from the erroneous Lemma 2 (ii) in \cite{Wer3} (see \cite{Wer13}), which, in particular, allows to deduce that $\Phi$ provides a construction for an equilibrium state for a contractive Markov system (see \cite{Wer10} and \cite{Wer12}) because it is absolutely continuous with respect to one (see Lemma 1 and Theorem 1 in \cite{Wer12}), where the existence of the latter is known through the Krylov-Bogolyubov argument. 

\begin{lemma}\label{dpfl}
	Let $\phi'_0$ and $\phi_0$ be non-negative measures on $\sigma$-algebra $\A_0$ such that $\phi'_0(X)<\infty$. Suppose $\phi'_0\ll\phi_0$  and $\phi_0\circ S^{-1}= \phi_0$. Then for the corresponding DDMs on $\B$ holds true the relation
	\[\Phi'\ll\Phi.\] 
\end{lemma} 
{\it Proof.}
Since $\phi_0\circ S^{-1}= \phi_0$, $\Phi|_{\A_m}=\phi_m$ for all $m\leq0$ by Proposition \ref{ip} (ii). Let $\epsilon>0$ and $\delta>0$ be such that $\phi'_0(A)<\epsilon/2$ whenever $\phi_0(A)<\delta$ for all $A\in\A_0$. Let $B\in\B$ such that $\Phi(B)<\delta$.  Then, by Proposition \ref{ip} (i), there exists  $(A_k)_{k\in\N}\subset\bigcup_{m\leq 0}\A_m$ such that $B\subset\bigcup_{k\in\N}A_k$ and $\sum_{k\in\N}\Phi(A_k)<\delta$. Then, by \eqref{ads}, for each $n\in\N$, there exists $m_n\leq 0$ such that $S^{m_n}\bigcup_{k=1}^nA_k\in\A_0$. Hence, by Theorem \ref{ddmt} (i), $\phi_0(S^{m_n}\bigcup_{k=1}^nA_k)=\Phi(\bigcup_{k=1}^nA_k)<\delta$ for all $n\in\N$, and therefore, $\phi'_0(S^{m_n}\bigcup_{k=1}^nA_k)<\epsilon/2$ for all $n\in\N$. Thus, by the $S$-invariance of $\Phi'$ on $\B$,
\begin{eqnarray*}
	\Phi'(B)&\leq&\lim\limits_{n\to\infty}\Phi'\left(S^{m_n}\bigcup_{k=1}^nA_k\right)\\
	&\leq&\limsup\limits_{n\to\infty}\phi'_0\left(S^{m_n}
	\bigcup\limits_{k=1}^nA_k\right) \\
	&<&\epsilon.
\end{eqnarray*}
\hfill$\Box$

The inference on the relation between $\Phi'|_\B$ and $\Phi|_\B$ from $\phi'_0\ll\phi_0$  if $\phi'_0\circ S^{-1}= \phi'_0$ is more subtle. To that is devoted another article \cite{Wer15}, which requires the measure theory developed here.

\section{Examples}\label{es}

Although it is easy to give an example of $\Phi(X)=0$ using an atomic $\phi_0$, the first example shows also that  
the atomicity of the initial measure does not imply $\Phi(X)=0$. It builds up on Example 1 in \cite{Wer10}.

\begin{Example}\label{ex1}
	Let $X := \{0,1\}^{\mathbb{Z}}$ (be the set of all $(...,\sigma_{-1},\sigma_0,\sigma_{1},...)$, $\sigma_i\in\{0,1\})$ and $S$ be the left shift map on $X$ (i.e. $(S\sigma)_i=\sigma_{i+1}$ for all $i\in\Z$). Let ${_0[a]}$ denote a cylinder set (i.e. the set of all $(\sigma_i)_{i\in\Z}\in X$ such that $\sigma_0=a$ where $a\in\{0,1\}$). Set $\A:=\{\emptyset,{_0[0]},{_0[1]},X\}$. Let $\sigma'\in X$ be given by
	\begin{equation*}
	\sigma'_i:=\left\{\begin{array}{cc} 0&  \mbox{if }i \mbox{ is even }\\
	1& \mbox{ otherwise }
	\end{array}\right.
	\end{equation*}
	for all $i\in\mathbb{Z}$.
	Let $\phi_0$ be the  measure on $\mathcal{A}_0$ given by
	\begin{equation*}
	\phi_0(A):=1_A(\sigma')\mbox{ for all }A\in\mathcal{A}_0.
	\end{equation*}
	
	Then $\Phi(X) = 0$, since $(...,\ \emptyset,\ \emptyset,  \ _0[0], \ _0[1])\in\C(X)$. Set
	\[\phi_{n0}:=\frac{1}{n+1}\sum\limits_{0\leq i\leq n}\phi_0\circ S^{-i}\mbox{ for }n\in\mathbb{N},\]
	and let $\Phi^{(n)}$ be the corresponding DDM.
	Then $\phi_{10}$ is shift-invariant and $\phi_{n0} = \phi_{10}$ for all odd $n$. So, $\Phi^{(n)}(X) = 1$ for all odd $n$. For every even $n$, $\phi_{n0}\geq n/(n+1)\phi_{10}$. Thus $\Phi^{(n)}(X) \geq n/(n+1)$ for all even $n$.
\end{Example}

A natural field of applications for the theory is, of course,  the theory of Markov processes, where the initial measure $\phi_0$ is usually available. The next example is just a scratch in that direction. 
\begin{Example}
	Let $A:=(a_{ij})_{1\leq i,j\leq N}$ be an irreducible stochastic $N\times N$-matrix. Then there exists a unique probability measure $\pi$ on all subsets of $\{1,...,N\}$ such that $\pi A=\pi$, and it has the property  $\pi\{i\}>0$ for all $1\leq i\leq N$. Let $\pi^{(0)}$ be any other probability measure on all subsets of $\{1,...,N\}$ such that $\pi^{(0)}\{i\}>0$ for all $1\leq i\leq N$. Define
	\begin{equation}\label{dbmc}
	\lambda_0:=\min_{1\leq i\leq N}\left\{\frac{\pi\{i\}}{\pi^{(0)}\{i\}}\right\}\mbox{ and }\alpha_0:=\max_{1\leq i\leq N}\left\{\frac{\pi\{i\}}{\pi^{(0)}\{i\}}\right\}.
	\end{equation}
	Then
	\[\lambda_0\pi^{(0)}\leq \pi\leq\alpha_0\pi^{(0)}.\]
	Let $X:=\{1,...,N\}^{\Z}$ and $S$ be the left shift map on $X$. Let $\A$ be the $\sigma$-algebra on $X$ generated by the cylinder sets $_0[a]$, $a\in\{1,...,N\}$. Let $\phi_0$ and $\phi_0^{(0)}$ be the probability measures on $\A_0$ given by
	\[\phi_0\left({_0[i_1,...,i_n]}\right):=\pi\{i_1\}a_{i_1i_2}...a_{i_{n-1}i_n}\]
	and
	\[\phi_0^{(0)}\left({_0[i_1,...,i_n]}\right):=\pi^{(0)}\{i_1\}a_{i_1i_2}...a_{i_{n-1}i_n}\]
	for all $_0[i_1,...,i_n]\subset\{1,...,n\}^{\Z}$ and $n\geq0$. Then, obviously,
	\[\lambda_0\phi_0^{(0)}\left({_0[i_1,...,i_n]}\right)\leq \phi_0\left({_0[i_1,...,i_n]}\right)\leq\alpha_0\phi_0^{(0)}\left({_0[i_1,...,i_n]}\right)\]
	for all $_0[i_1,...,i_n]\subset\{1,...,n\}^{\Z}$ and $n\geq0$.
	Let $\Phi$ and $\Phi^{(0)}$ denote the DDMs resulting from $\phi_0$ and $\phi_0^{(0)}$ respectively. Let $Q\subset X$ and $(A_m)_{m\leq 0}\in\C(Q)$. Then
	\[\lambda_0\Phi^{(0)}\left(Q\right)\leq\lambda_0\sum\limits_{m\leq 0}\phi_m^{(0)}\left(A_m\right)\leq\sum\limits_{m\leq 0}\lambda_0\phi_0^{(0)}\left(S^mA_m\right)\leq\sum\limits_{m\leq 0}\phi_m\left(A_m\right).\]
	Hence,
	\begin{equation}\label{ubmch}
	\lambda_0\Phi^{(0)}\left(Q\right)\leq\Phi\left(Q\right).
	\end{equation}
	Similarly, one sees that
	\begin{equation}\label{lbmch}
	\Phi\left(Q\right)\leq\alpha_0\Phi^{(0)}\left(Q\right).
	\end{equation}
	Since $\phi_0\circ S^{-1} = \phi_0$, $\Phi(X) = 1$ by Proposition \ref{nnp}, and therefore,
	\[\Phi^{(0)}(X)\geq\frac{1}{\alpha_0}.\]
	Furthermore, \eqref{ubmch} and \eqref{lbmch} imply that
	\begin{equation}\label{ucmc}
	\left|\Phi^{(0)}\left(Q\right)-\Phi\left(Q\right)\right|\leq\max\left\{(\alpha_0-1),\left(\frac{1}{\lambda_0}-1\right)\right\}\ \ \ \mbox{ for all }Q\subset X.
	\end{equation}
	For example, let
	\[\pi^{(0)}=\frac{1}{N}\sum\limits_{j=1}^N\delta_j.\] Then
	\[\Phi^{(0)}\left(Q\right)\geq\frac{1}{\alpha_0}\Phi\left(Q\right)=\frac{1}{N\max\limits_{1\leq i\leq N}\pi\{i\}}\Phi\left(Q\right).\]
	Thus
	\[\Phi^{(0)}\left(X\right)\geq\frac{1}{N}.\]
	For any other $\pi^{(0)}$, there exists $k_0\in\N$ such that for every $k\geq k_0$,
	\[\pi^{(k)}:=\frac{1}{k}\sum\limits_{k=0}^{k-1}\pi^{(0)}A^k\]
	satisfies $\pi^{(k)}\{i\}>0$ for all $i\in\{1,...,N\}$. If $A$ is aperiodic, then one can take $\pi^{(k)}:=\pi^{(0)}A^k$  with such property. For $k\geq k_0$, let $\lambda_k$ and $\alpha_k$ be defined as in \eqref{dbmc} with $\pi^{(k)}$ in place of $\pi^{(0)}$. Then, since, by the Ergodic Theorem, $\lambda_k\to 1$ and $\alpha_k\to 1$, it follows by \eqref{ucmc} that
	\[\lim\limits_{k\to\infty}\Phi^{(k)}\left(Q\right)=\Phi\left(Q\right)\ \ \ \mbox{ for all }Q\subset X.\]
\end{Example}

For a more general example arising from Markov processes, where the essential boundedness of the density function is not that obvious, see \cite{Wer13}.

\subsection*{Acknowledgements} 
I would like to thank the editors of the journals Fundamenta Mathematicae and Mathematical Proceedings of the Cambridge Philosophical Society for making me aware of the poor level of the general measure-theoretic education. This compelled to add more explanations to the article, which most likely improved it.

\end{document}